\documentclass[12pt]{amsart} 
\usepackage{amsfonts,amsmath,amsthm,amssymb,latexsym,amsthm,newlfont,enumerate,color}

\newif\ifslide

\theoremstyle{plain}

\ifdefined\spacer
\newtheorem{theorem}{Theorem}
 \else
\newtheorem{theorem}{Theorem}[section]
\fi

\newtheorem{corollary}[theorem]{Corollary}
\newtheorem{lemma}[theorem]{Lemma}
\newtheorem{claim}[theorem]{Claim}

\newtheorem*{theorem*}{Theorem}

\newtheorem{proposition}[theorem]{Proposition}
\newtheorem{definition-lemma}[theorem]{Definition-Lemma}
\newtheorem{question}[theorem]{Question}
\newtheorem{red-question}[theorem]{\textcolor{red}{Question}}

\newtheorem{conjecture}[theorem]{Conjecture}

\theoremstyle{definition}
\newtheorem{definition}[theorem]{Definition}
\newtheorem{remark}[theorem]{Remark}

\newtheorem{example}[theorem]{Example}

\DeclareMathOperator{\lcm}{lcm}
\DeclareMathOperator{\Bsl}{Bs}
\DeclareMathOperator{\Cl}{Cl}

\newcommand{\m}{\mathcal}
\renewcommand{\P}{\mathbb P}

\def\ideal#1.{I_{#1}}
\def\ring#1.{\mathcal {O}_{#1}}

\def\fring#1.{\hat{\mathcal {O}}_{#1}}
\def\proj#1.{\mathbb {P}(#1)}
\def\pr #1.{\mathbb {P}^{#1}}
\def\dpr #1.{\hat{\mathbb {P}}^{#1}}
\def\af #1.{\mathbb A^{#1}}
\def\Hz #1.{\mathbb F_{#1}}
\def\Hbz #1.{\overline{\mathbb F}_{#1}}
\def\fb#1.{\underset #1 {\times}}
\def\rest#1.{\underset {\ \ring #1.} \to \otimes}
\def\au#1.{\operatorname {Aut}\,(#1)}
\def\deg#1.{\operatorname {deg } (#1)}

\def\pic#1.{\operatorname {Pic}\,(#1)}
\def\pico#1.{\operatorname{Pic}^0(#1)}
\def\picg#1.{\operatorname {Pic}^G(#1)}
\def\ner#1.{NS (#1)}
\def\rdown#1.{\llcorner#1\lrcorner}
\def\rfdown#1.{\lfloor{#1}\rfloor}
\def\rup#1.{\ulcorner{#1}\urcorner}
\def\rcup#1.{\lceil{#1}\rceil}

\def\n1#1.{\operatorname {N_1}(#1)}  
\def\cn1#1.{\overline{\operatorname {N^1}(#1)}} 
\def\cone#1.{\operatorname {NE}(#1)}     
\def\ccone#1.{\overline{\operatorname {NE}}(#1)}
\def\none#1.{\operatorname {NF}(#1)}
\def\cnone#1.{\overline{\operatorname {NF}}(#1)}
\def\mone#1.{\operatorname {NM}(#1)} 
\def\cmone#1.{\overline{\operatorname {NM}}(#1)}

\def\coef#1.{\frac{(#1-1)}{#1}}
\def\vit#1.{D_{\langle #1 \rangle}}
\def\mm#1.{\overline {M}_{0,#1}}
\def\H1#1.{H^1(#1,{\ring #1.})}
\def\ac#1.{\overline {\mathbb F}_{#1}}

\def\adj#1.{\frac {#1-1}{#1}}
\def\spn#1.{\overline{#1}}
\def\pek#1.#2.{\Cal P^{#1}(#2)}
\def\plk#1.#2.{\Cal P^{\leq #1}(#2)}
\def\ev#1.{\operatorname{ev_{#1}}}
\def\ilist#1.{{#1}_1,{#1}_2,\dots}
\def\bminv#1.{(\nu_1,s_1;\nu_2,s_2;\dots ;\nu_{#1},s_{#1};\nu_{r+1})}
\def\zinv#1.{(\nu_1,s_1;\nu_2,s_2;\dots ;\nu_{#1},s_{#1};0)}
\def\iinv#1.{(\nu_1,s_1;\nu_2,s_2;\dots ;\nu_{#1},s_{#1};\infty)}

\def\scr #1.{\mathcal #1}

\def\llist#1.#2.{{#1}_1,{#1}_2,\dots,{#1}_{#2}}
\def\ulist#1.#2.{{#1}^1,{#1}^2,\dots,{#1}^{#2}}
\def\lomitlist#1.#2.{{#1}_1,{#1}_2,\dots,\hat {{#1}_i}, \dots, {#1}_{#2}}
\def\lomitlistz#1.#2.{{#1}_0,{#1}_1,\dots,\hat {{#1}_i}, \dots, {#1}_{#2}}
\def\loc#1.#2.{\Cal O_{#1,#2}}
\def\fderiv#1.#2.{\frac {\partial #1}{\partial #2}}
\def\deriv#1.#2.{\frac {d #1}{d #2}}

\def\map#1.#2.{#1 \longrightarrow #2}
\def\rmap#1.#2.{#1 \dasharrow #2}
\def\emb#1.#2.{#1 \hookrightarrow #2}
\def\non#1.#2.{\text {Spec }#1[\epsilon]/(\epsilon)^{#2}}
\def\Hi#1.#2.{\text {Hilb}^{#1}(#2)}
\def\sym#1.#2.{\operatorname {Sym}^{#1}(#2)}
\def\Hb#1.#2.{\text {Hilb}_{#1}(#2)}
\def\Hm#1.#2.{\Hom_{#1}(#2)}
\def\prd#1.#2.{{#1}_1\cdot {#1}_2\cdots {#1}_{#2}}
\def\Bl #1.#2.{\operatorname {Bl}_{#1}#2}
\def\pl #1.#2.{#1^{\otimes #2}}
\def\mgn#1.#2.{\overline {M}_{#1,#2}}
\def\ialist#1.#2.{{#1}_1 #2 {#1}_2, #2\dots}
\def\pair#1.#2.{\langle #1, #2\rangle}
\def\vandermonde#1.#2.{\left|
\begin{matrix}
1 & 1 & 1 & \dots & 1\\
{#1}_1 & {#1}_2 & {#1}_3 & \dots & {#1}_{#2}\\
{#1}_1^2 & {#1}_2^2 & {#1}_3^2 & \dots & {#1}_{#2}^2\\
\vdots & \vdots & \vdots & \ddots & \vdots\\
{#1}_1^{#2-1} & {#1}_2^{#2-1} & {#1}_2^{#2-1} & \dots & {#1}_{#2}^{#2-1}\\
\end{matrix}
\right|
}
\def\vandermondet#1.#2.{\left|
\begin{matrix}
1 & {#1}_1   & {#1}_1^2 & \dots & {#1}_1^{#2-1}\\
1 & {#1}_2   & {#1}_2^2 & \dots & {#1}_2^{#2-1}\\
1 & {#1}_3   & {#1}_3^2 & \dots & {#1}_3^{#2-1}\\
\vdots & \vdots & \vdots & \ddots & \vdots\\
1 & {#1}_{#2}& {#1}_{#2}^2 & \dots & {#1}_{#2}^{#2-1}\\
\end{matrix}
\right|
}
\def\gr#1.#2.{\mathbb{G}(#1,#2)}

\def\alist#1.#2.#3.{{#1}_1 #2 {#1}_2 #2\dots #2 {#1}_{#3}}
\def\zlist#1.#2.#3.{#1_0 #2 #1_1 #2\dots #2 #1_{#3}}
\def\lomitlist30#1.#2.#3.{{#1}_0,{#1}_1 #2 \dots #2\hat {{#1}_i} #2\dots #2 {#1}_{#3}}
\def\lmap#1.#2.#3.{#1 \overset{#2}{\longrightarrow} #3}
\def\mes#1.#2.#3.{#1 \longrightarrow #2 \longrightarrow #3}
\def\ses#1.#2.#3.{0\longrightarrow #1 \longrightarrow #2 \longrightarrow #3 \longrightarrow 0}
\def\les#1.#2.#3.{0\longrightarrow #1 \longrightarrow #2 \longrightarrow #3}
\def\res#1.#2.#3.{#1 \longrightarrow #2 \longrightarrow #3\longrightarrow 0}
\def\Hi#1.#2.#3.{\text {Hilb}^{#1}_{#2}(#3)}
\def\ten#1.#2.#3.{#1\underset {#2}{\otimes} #3}
\def\lomitlist30#1.#2.#3.{{#1}_0 #2 {#1}_1 #2 \dots #2 \hat {{#1}_i} #2 \dots #2 {#1}_{#3}}
\def\mderiv#1.#2.#3.{\frac {d^{#3} #1}{d #2^{#3}}}

\def\Hom{\operatorname{Hom}}

\def\dim{\operatorname{dim}}

\def\deg{\operatorname{deg}}

\def\Pic{\operatorname{Pic}}

\def\Sing{\operatorname{Sing}}

\def\rk{\operatorname{rk}}

\def\rest{\operatorname{res}}

\def\C{\mathbb C}

\def\e{\Cal E}

\def\e1{E_1}
\def\e2{E_2}

\def\Z{\mathbb Z}
\def\N{\mathbb N}

\def\mapdown#1{\big\downarrow\rlap{$\vcenter{\hbox{$\scriptstyle#1$}}$}}

\def\mapse#1{
{\vcenter{\hbox{$\mathop{\smash{\raise1pt\hbox{$\diagdown$}\!\lower7pt
\hbox{$\searrow$}}\vphantom{p}}\limits_{#1}\vphantom{\mapdown{}}$}}}}

\def\VR#1.{height#1pt&\omit&&\omit&&\omit&&\omit&&\omit&\cr}

\def\VRT#1.{height#1pt&\omit&&\omit&\cr}

\usepackage[colorlinks=true,pagebackref,hyperindex]{hyperref}
\usepackage{amsmath}
\usepackage{amssymb}
\usepackage{amsthm}
\usepackage[all,cmtip]{xy}

\title{Effective non-vanishing for Fano weighted complete intersections}
\date{\today}

\author{Marco Pizzato}
\address{}
\email{marco.pizzato1@gmail.com}

\author{Taro Sano}
\address{Department of Mathematics, Graduate School of Science, 
Kobe university, 
1-1, Rokkodai, Nada-ku, Kobe 657-8501, Japan} 
\email{tarosano@math.kobe-u.ac.jp} 

\author{Luca Tasin}
\address{Universit\`a Roma Tre, Dipartimento di Matematica e Fisica, Largo San Leonardo Murialdo I-00146 Roma, Italy} 
\email{ltasin@mat.uniroma3.it}

\begin{document}



\begin{abstract}
We show that Ambro--Kawamata's non-vanishing conjecture holds true for a quasi-smooth WCI $X$ which is Fano or Calabi-Yau, i.e.\ we prove that, if $H$ is an ample Cartier divisor on $X$, then $|H|$ is not empty. If $X$ is smooth, we further show that the general element of $|H|$ is smooth. 
We then verify Ambro--Kawamata's conjecture for any quasi-smooth weighted hypersurface.  
We also verify Fujita's freeness conjecture for a Gorenstein quasi-smooth weighted hypersurface. 

For the proofs, we introduce the arithmetic notion of regular pairs and enlighten some interesting connection with the Frobenius coin problem.
\end{abstract}

\maketitle



\section{Introduction}

Complete intersections in weighted projective spaces (WCI's for short) form a natural class of varieties which are particularly interesting from the point of view of higher dimensional algebraic geometry. We refer to \cite{Dolgachev}, \cite{Mori75} and \cite{Dimca86} for a general treatment of these varieties.




Reid \cite{Reid80,Reid87} and Iano-Fletcher \cite{Fletcher00} systematically investigated notable examples of WCI's and started their classification.  Several results have been afterwards obtained concerning boundedness and classification, see for example \cite{Johnson-Kollar}, \cite{Chen:2012aa}, \cite{Chen:2009aa} and \cite{PS16}. 





\smallskip

The main motivation of this paper is to study the following conjecture in the realm of WCI's, in particular for what concerns the case of Fano and Calabi-Yau varieties.

\begin{conjecture}[Ambro--Kawamata]\label{conj:Ambro-Kawamata}
Let $(X,\Delta)$ be a klt pair and $H$ be an ample Cartier divisor on $X$ such that $H-K_X-\Delta$ is ample.
Then $|H| \ne \emptyset$.
\end{conjecture}

For an introduction to this conjecture, see \cite{Ambro99} and in particular \cite{Kawamata00} where the 2-dimensional case is proven. 
In the smooth setting, Ionescu \cite[Page 321]{Cetraro93} and Beltrametti--Sommese \cite{BS95} proposed related conjectures.

Ambro-Kawamata's conjecture is known to be true in full generality only in dimension 1 and 2. Several cases have been studied, especially in dimension 3 (see for instance \cite{Xie09},	\cite{BH10}, \cite{Horing12} and \cite{CaoJiang}).

A {\it fundamental divisor} on a Fano variety $X$ is an ample Cartier divisor $H$ which is primitive 
and proportional to $-K_X$. 
In the classification of Fano varieties, it is important to investigate the properties of the general member of the linear system given by $H$, see for instance \cite{Ambro99}.  
The second purpose of this note is to study this problem in the case of Fano and Calabi-Yau smooth WCI's.  



\vspace{2mm}

The main result of this paper is the following.

\begin{theorem}\label{thm:non-vanishing}
Let $X=X_{d_1,\ldots,d_c} \subset \P(a_0,\ldots,a_n)$ be a well-formed quasi-smooth weighted complete intersection which is not a linear cone and $H$ be an ample Cartier divisor on $X$. Assume that $X$ is Fano or Calabi-Yau. Then $|H| \neq \emptyset$. 

Moreover, if $X$ is smooth, then the number of $a_i=1$ is at least $c$ and the general element of $|\mathcal O_X(1)|$ is smooth.	
\end{theorem}

For a smooth Fano WCI, it was already proved in \cite[Lemma 3.3]{PS16} that at least two weights 
are $1$, which implies the non-vanishing for a smooth Fano WCI.  
In addition, it is easy to prove Conjecture \ref{conj:Ambro-Kawamata} for any smooth WCI of codimension 1 and 2, see Remark \ref{rmk:codimension2}.

It is particularly interesting that, in the smooth case, 
we can prove the smoothness of the general member of the fundamental linear system (Corollary \ref{cor:smooth} (ii)). 

Theorem \ref{thm:non-vanishing} is a direct consequence of Corollaries \ref{cor:smooth} and \ref{cor:singular}. 
In particular, in Corollary \ref{cor:smooth}, we show that if $X=X_{d_1,\ldots,d_c} \subset \P(a_0,\ldots,a_n)$ is a smooth well-formed Fano WCI which is not a linear cone, then the number of $i$ for which $a_i=1$ is at least $c+1$.  
By using this, we can then show that the general element of $|\m O_X(1)|$ is quasi-smooth (from which smoothness follows easily). 
One can not expect a similar statement for a general member of the fundamental linear system of a singular quasi-smooth WCI, as Example \ref{ex:funddiv} shows. 
We also give a description of the base locus of $|\m O_X(1)|$ in Remark \ref{rem:baselocus} and an example 
whose base locus $\Bsl |\m O_X(1)|$ is singular and not quasi-smooth in Example \ref{ex:baselocus}. 

In \cite[Corollary 4.2]{Przyjalkowski:2017ab}, the authors show that for a smooth well-formed Fano WCI $X$ the number of $a_i=1$ is at least $I(X)=\sum a_i -\sum d_j$ when $c \le 2$ and write that they expect this to hold for any codimension. As a consequence of Proposition \ref{prop:regular}, we can confirm this expectation, see Corollary \ref{cor:I(X)}.

\smallskip

In the case of hypersurfaces, we can prove the following stronger result, which is the combination of Propositions \ref{prop:hypnonvanishing} and \ref{prop:Fujita}.

\begin{theorem}\label{thm:hypersurfaces}
Let $X=X_d \subset \P=\P(a_0,\ldots,a_n)$ be a well-formed quasi-smooth hypersurface of degree $d$ which is not a linear cone. 
\begin{enumerate}
	\item If $H$ is an ample Cartier divisor on $X$ such that $H - K_X$ is ample, then $|H|$ is not empty.
  \item If $X$ is Gorenstein and $H$ is an ample Cartier divisor, then $K_X+mH$ is globally generated for any $m \ge n$.
\end{enumerate}
\end{theorem}

The second part of the statement is known as Fujita's freeness conjecture and it has been proven in the smooth setting up to dimension 5 (see \cite{Reider88}, \cite{EL93}, \cite{Kawamata97} and \cite{YZ15}). 

\smallskip

\subsection{The methods}
The above theorems are obtained by studying the arithmetic properties of quasi-smooth WCI's.  
More precisely, in Section \ref{sec:section0}, we prove a criterion (Proposition \ref{prop:qsmcrit}) for a WCI to be quasi-smooth, which generalizes Iano-Fletcher's criterion in codimension 1 and 2 (see \cite[Sect. 8]{Fletcher00}). 
We then exploit some arithmetic consequences of quasi-smoothness. In particular, Proposition \ref{prop:i1i2} motivates the introduction of an {\it $h$-regular pair} (see Definition \ref{def:h-regular}) which turns out to be a key tool in our treatment. 

Given a positive integer $h$, a pair $(d;a)=(d_1,\ldots,d_c; a_0,\ldots,a_n) \in \N^{c} \times \N^{n+1}$ is called $h$-regular if for any $I= \{i_1, \ldots,  i_k \} \subset \{0,\ldots,n\}$ such that $a_I:=\gcd(a_{i_1}, \ldots, a_{i_k})>1$, either $a_I \mid h$ or there are  distinct integers $p_1, \ldots, p_k$ such that 
$$ 
a_I \mid d_{p_1}, \ldots, d_{p_k}.
$$    

Set $\delta(d;a):=\sum_{j=1}^c d_j - \sum_{i=0}^n a_i$. By Proposition \ref{prop:i1i2}, any quasi-smooth (well-formed) WCI $X=X_{d_1,\ldots,d_c} \subset \P(a_0,\ldots,a_n)$ gives rise to an $h$-regular pair $(d;a)=(d_1,\ldots,d_c;a_0,\ldots,a_n)$, where $h$ is the smallest positive integer for which $\m O_X(h)$ is Cartier.  Remembering that $K_X= \m O_X(\delta)$, the non-vanishing for a Fano or Calabi-Yau WCI follows from Proposition \ref{prop:h-regular}, which says that, if $(d;a)$ is $h$-regular such that $a_i \neq d_j$ and $a_i \nmid h$ for any $i,j$, then $\delta(d;a) >0$. A more accurate statement (Corollary \ref{cor:smooth}) is needed to prove that, if $X$ is smooth, then the general element of $|\m O_X(1)|$ is also smooth.

We now spend some words for the case $h=1$. In this case, the pair $(d;a)$ is simply called \emph{regular}. 
A smooth WCI $X$ gives rise to a regular pair $(d;a)$. The non-vanishing is then equivalent to prove that 
$$
\delta(d; a) \ge G(a_0,\ldots,a_n),
$$
where $G(a_0,\ldots,a_n)$ is the Frobenius number of $a_0,\ldots,a_n$, i.e.\ the greatest integer which is not a non-negative integral combination of $a_0,\ldots,a_n$. In Conjecture \ref{conj:regular}, we speculate that $\delta (d; a) \ge G(a_0,\ldots,a_n)$ for a regular pair $(d;a)$, under some natural assumptions. This would imply Ambro--Kawamata's conjecture for any smooth WCI. 

We believe that this conjecture is interesting also from the arithmetic point of view, since it would give new bounds for the Frobenius number (see Section \ref{sec:Frobenius} for details).


\section{Preliminaries and notation}

In this section, we recall some basic facts about weighted complete intersections and fix our notation. See \cite{Dolgachev} or \cite{Fletcher00} for further details.

Let $\mathbb{N}$ (resp.\ $\mathbb{N}_+$) be the set of non-negative (resp.\ positive) integers. 
Let $a_0,\ldots,a_n \in \mathbb{N}_+$. 
We define $\P:=\P(a_0,\ldots,a_n)$ to be the weighted projective space with weights $a_0, \ldots, a_n$, i.e.\ $\P= \mathrm{Proj}\ \C[x_0,\ldots,x_n]$, where $x_i$ has weight $a_i$. 
We denote $$\mathbb{P}(\underbrace{b_1, \ldots, b_1}_{k_1}, \ldots, \underbrace{b_l, \ldots, b_l}_{k_l})$$ by $\mathbb{P}(b_1^{(k_1)}, \ldots, b_l^{(k_l)})$ for short.

Note that if we start with $x_0, \ldots, x_n$ to be affine coordinates on $\mathbb A^{n+1}$ and $\C^*$ acting on $\mathbb A^{n+1}$ via
$$
\lambda \cdot (x_0, \ldots, x_n)= (\lambda^{a_0}x_0, \ldots, \lambda^{a_n}x_n)
$$
for any $\lambda \in \C^*$, then $\P$ is just the quotient $(\mathbb A^{n+1} \setminus \{0\})/\C^*$.

We always assume that $\P$ is \emph{well-formed}, i.e.\  the greatest common divisor of any $n$ weights is 1. 
For any $I= \{i_1, \ldots, i_k \} \subset \{0,\ldots,n\}$, the stratum $\Pi_I$ is defined as 
$$
\Pi_I:= \{ x_i=0 : i \notin I\}.
$$
The singular locus of $\P$ is the union of all strata $\Pi_I$ for which $a_I:=\gcd(a_i)_{i \in I} >1$. Any point of the interior $\Pi^0_I$ of a stratum $\Pi_I$ is locally isomorphic to a quotient singularity of type 
$$
\frac{1}{a_I}(a_0, \ldots, \hat{a_{i_1}}, \ldots, \hat{a_{i_k}}, \ldots,  a_n) \times \mathbb{C}^{k-1}. 
$$ 
Here, for $r \in \mathbb{N}_+$ and $a_1, \ldots, a_n \in \mathbb{N}$ such that 
$ \gcd (r, a_1, \ldots, a_n)=1$, a {\it quotient singularity} of type $1/r (a_1, \ldots, a_n)$ means a quotient $\mathbb{C}^n/ \mathbb{Z}_r$ by the action of 
a cyclic group $\mathbb{Z}_r$ of order $r$ as $g \cdot z_i = \zeta_r^{a_i} z_i$ for $i=1, \ldots , n$, 
where $g \in \mathbb{Z}_r$ is a generator and $\zeta_r$ is an $r$-th primitive root of unity. 
We also denote by $\mathbb{C}^n/ \mathbb{Z}_r(a_1, \ldots, a_r)$ this quotient affine variety. 
Let $\pi \colon \mathbb{C}^n \rightarrow U:= \mathbb{C}^n/ \mathbb{Z}_r(a_1, \ldots, a_r)$ 
be the quotient morphism. We have an eigen-decomposition 
\[
\pi_* \mathcal{O}_{\mathbb{C}^n} = \bigoplus_{i=0}^{r-1} \mathcal{F}_i, 
\]
where $\mathcal{F}_i := \{  f \in \mathcal{O}_{\mathbb{C}^n} \mid g \cdot f = \zeta_r^i f \}$ is 
the $\mathcal{O}_U$-submodule of $\pi_* \mathcal{O}_{\mathbb{C}^n}$ consisting of 
$\mathbb{Z}_r$-eigenfunctions of eigenvalue $\zeta_r^i$. 
Note that $\mathcal{F}_i \simeq \mathcal{O}_U(D_f)$, where a divisor $D_f:= (f=0)/ \mathbb{Z}_r$ on $U$ is 
defined by a function $f \in \mathcal{F}_i$.

\begin{proposition}
The divisor class group of $U:= \mathbb{C}^n/ \mathbb{Z}_r(a_1, \ldots, a_r)$ is 
\begin{equation}\label{eq:classgroup}
\Cl  U \simeq \mathbb{Z}_r \cdot \mathcal{F}_1. 
\end{equation}
 \end{proposition}

\begin{proof} 
We have an inclusion $\iota \colon \mathbb{Z}_r \cdot \mathcal{F}_1 \hookrightarrow \Cl U$. 
It is enough to show that this is surjective. 
Let $D \subset U$ be a prime divisor. Then $\pi^{-1}(D)$ is a divisor on $\mathbb{C}^n$ 
defined by some $\mathbb{Z}_r$-eigenfunction $f_D \in \mathcal{O}_{\mathbb{C}^n}$. 
Let $i(D) \in \mathbb{Z}_r$ be an element such that $g \cdot f_D = \zeta_r^{i(D)} f$. 
Then we can check that $\mathcal{O}_U(D) \simeq \mathcal{F}_{i(D)}$. 
Since $g^{i} \cdot \mathcal{F}_1 \simeq \mathcal{F}_i$ for $i \in \mathbb{Z}_r$, 
we see that $\iota$ is surjective.   

We can also check the isomorphism by toric computation. 
Since $\mathbb{C}^n/ \mathbb{Z}_r(a_1, \ldots, a_n)$ is a toric variety, 
 we can compute its class group by using the information of the cone and lattice. 
 (cf.\ \cite[3.4 Proposition]{Fulton93}) 
 More precisely, it is the quotient $\mathbb{Z}^n /M$, where 
 $M:= \{(m_1, \ldots, m_n) \in \mathbb{Z}^n \mid \sum_{i=1}^n m_i a_i \equiv 0 \mod r \}$. 
\end{proof}

\begin{definition}
Let $X$ be a (closed) subvariety of codimension $c$ in $\P$. 
Then $X$ is  \emph{well-formed} if 
$$
\mathrm{codim}_X (X \cap \mathrm{Sing}(\P)) \ge 2.
$$
Let $\pi: \mathbb{A}^{n+1}\setminus \{0\} \to \P$ be the natural projection. Then $X$ is \emph{quasi-smooth} if $\pi^{-1}(X)$ is smooth.

The variety $X$ is said to be a \emph{weighted complete intersection} (WCI for short) of multidegree $(d_1,\ldots,d_c)$ if its weighted homogenous ideal in $\mathbb C[x_0,\ldots,x_n]$ is generated by a regular sequence of homogenous polynomials $\{f_j\}$ such that $\deg f_j=d_j$ for $j=1,\ldots,c$.
We denote by $X_{d_1,\ldots,d_c}$ a general element of the family of WCI's of multidegree $(d_1,\ldots,d_c)$.

Finally, $X_{d_1,\ldots,d_c} \subset \P$ is said to be a \emph{linear cone} if $d_j=a_i$ for some $i$ and $j$.
\end{definition}

Note that by \cite[Prop. 8]{Dimca86}, if $X$ is a well-formed quasi-smooth WCI, then 
$$
\mathrm{Sing}(X)=X \cap \mathrm{Sing}(\P).
$$

\begin{proposition}
If $X$ is a quasi-smooth WCI of dimension $\ge 3$ , then its divisor class group is a free  $\Z$-module generated by $\m O_X(1)$, where $\m O_X(1):=\m O_{\P}(1)_{| X}$. (We freely mix the divisorial and the sheaf notation. )
\end{proposition}

\begin{proof}
The proof is same as \cite[Lemma 3.5]{MR1798983}. 
This follows from the parafactoriality of an l.c.i.\ local ring (\cite{MR1266175}). 
\end{proof}

If $X \subset \P$ is a well-formed quasi-smooth WCI, then 
$$
\omega_X=K_X =\m O_X(\sum_{j=1}^c d_j - \sum_{i=0}^n a_i),
$$ 
see \cite[Thm. 3.3.4]{Dolgachev}. We usually write $\delta :=\sum_{j=1}^c d_j - \sum_{i=0}^n a_i$.

\smallskip

	The following result shows that the dimension of the linear system $|\mathcal{O}_X(n)|$ can be computed by the weights of the coordinates. 
	
	\begin{lemma}(\cite[Lemma 7.1]{Fletcher00})\label{lem:sectiondim} 
		Let $X \subset \mathbb{P}(a_0, \ldots, a_n)$ be a well-formed quasi-smooth WCI. 
		Let $A := k [x_0, \ldots, x_n]/(f_1, \ldots, f_c) $ be the homogeneous coordinate ring of $X$ and 
		$A_k$ be the $k$-th graded part for $k \in \mathbb{Z}$. 
		
		Then we have 
		\[
		H^0(X, \mathcal{O}_X(k)) \simeq A_k. 
		\]  
		\end{lemma}
	
	\begin{proof}
		See also \cite[3.4.3]{Dolgachev}. This follows since the homogeneous coordinate ring $A$
		  is Cohen-Macaulay 
		and we have $H^1_{\mathfrak{m}}(A) =0$, where $\mathfrak{m}:=(x_0, \ldots, x_n)$ is the maximal ideal.  
		\end{proof}


\section{Properties of quasi-smooth WCI's} \label{sec:section0}


In the following proposition, we give a necessary and sufficient condition for quasi-smoothness 
of a WCI. 

\begin{proposition}\label{prop:qsmcrit} 
	Let $X= X_{d_1, \ldots, d_c} \subset \mathbb{P}(a_0, \ldots, a_n)=: \mathbb{P}$ be a
	quasi-smooth WCI which is not a linear cone. 
	Let $x_0, \ldots, x_n$ be the coordinates of 
	$\mathbb{P}(a_0, \ldots, a_n)$.  
Fix $I:= \{i_1, \ldots, i_k \} \subset \{0,\ldots,n\}$ and let $\rho_I:= \min \{c, k \}$. 
For $m=(m_1, \ldots, m_k)$, let $x_I^{m}:= \prod_{j=1}^k x_{i_j}^{m_j}$. 
For a finite set $A$, let $|A|$ be the number of its elements. 
Then one of the following holds. 

\begin{enumerate}
	\item[(Q1)] There exist distinct integers $p_1,\ldots, p_{\rho_I} \in \{1, \ldots, c \}$ 
	and $M_1, \ldots, M_{\rho_I} \in \mathbb{N}^{k}$ 
	such that the monomial $x_I^{M_{j}}$ has the degree $d_{p_j}$ 
	for $j=1, \ldots, \rho_I$. 
	
	\item[(Q2)] There exist a permutation $p_1, \ldots, p_c$ of $\{1, \ldots, c \}$, an integer $l <  \rho_I$, and 
	integers $e_{\mu, j} \in \{0, \ldots, n \} \setminus I$ for $\mu=1, \ldots, k-l$ and $j= l+1, \ldots, c$ 
	such that there are monomials $x_I^{M_j}$ of degree $d_{p_j}$ 
	for $j=1, \ldots, l$ and distinct $k-l$ monomials 
	$\{ x_{e_{\mu, j}} x_I^{M_{\mu, j}} : \mu= 1, \ldots, k-l \}$ of degree $d_{p_j}$ for each $j=l+1, \ldots, c$  
	which satisfy the following: for any subset $J \subset \{l+1, \ldots, c \}$, we have $\left| \{e_{\mu,j} : j \in J, \mu = 1, \ldots, k-l \} \right|  \ge k-l + |J| -1$.  
	\end{enumerate} 
	Conversely, if we have (Q1) or (Q2) for all $I$, then a general WCI 
	$X_{d_1, \ldots, d_c} \subset \mathbb{P}(a_0, \ldots, a_n)$ is quasi-smooth. 
	\end{proposition}

\begin{remark}
This generalizes \cite[Theorem 8.7]{Fletcher00} in codimension 2 case. 
A weaker necessary condition for the quasi-smoothness is written in \cite[Proposition 2.3]{Chen:2012aa}. 
Although we shall not use the new part of Proposition \ref{prop:qsmcrit} in the main part of this paper, 
we believe it is an interesting result on its own. 
\end{remark}

\begin{proof} The framework of the proof is similar to that of \cite[Theorem 8.7]{Fletcher00}. 
	
	Let $F_j := |\mathcal{O}_{\mathbb{P}}(d_j)|$ be the linear system of weighted homogeneous polynomials of degree $d_j$. 
	For $j=1, \ldots, c$, let $f_j$ be a general homogeneous polynomial of degree $d_j$ such that 
	$X= (f_1=\cdots = f_c =0) \subset \mathbb{P}(a_0, \ldots, a_n)$. 
	Let $C^*_X \subset \mathbb{A}^{n+1} \setminus \{0 \}$ be the cone over $X$ defined by the polynomials 
	$f_1, \ldots, f_c$ with the following diagram 
	\[
	\xymatrix{
	C_X^* \ar@{^{(}->}[r] \ar[d] & \mathbb{A}^{n+1} \setminus \{0 \} \ar[d] \\
	X \ar@{^{(}->}[r] & \mathbb{P}. 
	}
	\]
	Without loss of generality, we may assume $I = \{ 0, \ldots, k-1 \}$ in the statement. 
	Let $\Pi:= (x_k = \cdots = x_n =0) \subset \mathbb{A}^{n+1}$ be the stratum corresponding to $I$ 
	and $\Pi^0 \subset \Pi $ be the open toric stratum. 
	By expanding $f_{\lambda}$ for $\lambda =1, \ldots, c$ in terms of $x_k, \ldots, x_n$, we can write 
	\[
	f_{\lambda} = h_{\lambda}(x_0, \ldots, x_{k-1}) + \sum_{i=k}^n x_i g^i_{\lambda}(x_0, \ldots, x_{k-1}) + R_{\lambda}(x_0, \ldots, x_n),  
	\]
	where $h_{\lambda}, g^i_{\lambda} \in \mathbb{C}[x_0, \ldots, x_{k-1}]$ and $R_\lambda \in \mathbb{C}[x_0, \ldots, x_n]$ satisfies $\deg_{x_k, \ldots, x_n} R_{\lambda} \ge 2$. 
	
	Note that $X$ is quasi-smooth if and only if $C_{X}^*$ is smooth along all the coordinate strata. 
	We shall show that $C_X^*$ is smooth along $\Pi^0$ when either (Q1) or (Q2) holds for $I$. 
	Let $\rho:= \rho_I$ for short. 

	Suppose that (Q1) holds. Then $h_{p_1}, \ldots, h_{p_{\rho}}$ are nonzero on $\Pi^0$. 
	If some of $h_{p_j}$ involves only one monomial, then we have $\Pi^0 \cap C_{X}^* = \emptyset$. 
	So we may assume that each of $h_{p_1}, \ldots, h_{p_{\rho}}$ involves at least $2$ monomials. 
	Thus we see that the linear systems $F_{d_{p_1}}, \ldots, F_{d_{p_{\rho}}}$ do not have base locus on $\Pi^0$. 
	By Bertini's theorem, we see that $(f_{p_1} = \cdots = f_{p_{\rho}} =0) \subset \mathbb{A}^{n+1}$  
	is smooth along $\Pi^0$ when $k \ge c$. 
	When $k < c$, we have $(f_{p_1} = \cdots = f_{p_{\rho}} =0) \cap \Pi^0 = \emptyset$. 
	Therefore $C_X^*$ is nonsingular along $\Pi^0$. 
	 
	Next suppose that (Q2) holds. 
	By permutation, we may assume that $p_i = i$. 
	Then $h_1, \ldots, h_{l}$ are nonzero on $\Pi^0$.  
	Hence the base locus of $F_{d_\lambda}$ is disjoint with $\Pi^0$  for  $\lambda = 1, \ldots, l$. 
	By Bertini's theorem, we see that $(f_1 = \cdots = f_l =0)$ is nonsingular along $\Pi^0$. 
	 We may assume that the Jacobian of $(f_1 = \cdots = f_c =0) \subset \mathbb{A}^{n+1}$ at $P \in \Pi^0$ 
	 is of the form 
	 \[
	 \begin{pmatrix}
	 \frac{\partial f_1}{\partial x_0} & \cdots & \frac{\partial f_1}{\partial x_{k-1}} & & & \\ 
	 \vdots & & \vdots & & \text{{\huge *}} & \\ 
	 \frac{\partial f_l}{\partial x_0} & \cdots & \frac{\partial f_l}{\partial x_{k-1}} & & & \\ 
	  & & & g^{k}_{l+1} & \cdots & g^n_{l+1} \\ 
	  & \text{\huge 0} & & \vdots & & \vdots \\ 
	   & & & g^{k}_{c} & \cdots & g^n_{c}
	 \end{pmatrix} (P) 
	 \] 
	 since we 
	 have   $h_{\lambda} =0$ for $\lambda = l+1, \ldots, c$. 
	 Note that the block matrix 
	 \[
	 \begin{pmatrix}
	 \frac{\partial f_1}{\partial x_0} & \cdots & \frac{\partial f_1}{\partial x_{k-1}} \\ 
	 \vdots & & \vdots \\ 
	 \frac{\partial f_l}{\partial x_0} & \cdots & \frac{\partial f_l}{\partial x_{k-1}} \\
	 \end{pmatrix}(P)
	 \] 
	 has maximal rank $l$ at $P \in \Pi^0$ since $(f_1 = \cdots = f_l =0)$ is nonsingular along $\Pi^0$. 
	 Hence it is enough to show that the matrix
	 \begin{equation}\label{eq:M_P} 
	 M_P:= 
	 \begin{pmatrix}
	 g_{l+1}^k & \cdots & g_{l+1}^n \\ 
	 \vdots & & \vdots \\ 
	 g^k_c & \cdots & g_c^n 
	 \end{pmatrix} (P) 
	 \end{equation}
	 has maximal rank $c-l$. 
	 
	 Note that  there are at least $k-l$  elements of $K_{\lambda}:= \{ i \in \{k, \ldots, n \} :  g^i_{\lambda} \neq 0 \}$ for each $\lambda = l+1, \ldots, c$. 
	 By $ |K_{\lambda}| \ge k-l$, we see that each row vector of $M_P$ is nonzero for $P \in \Pi^0$. 
	 Indeed, for each $\lambda' =l+1, \ldots, c$, the intersection 
	 \[
	  \bigcap_{\lambda =1}^l (h_{\lambda} =0) \cap \bigcap_{i=k}^n (g^i_{\lambda'} = 0) \cap \Pi^0 
	 \]
	  is contained in at least $k=l+(k-l)$ free linear systems on $k$-dimensional $\Pi^0$, and it is empty. 
	 Thus we may assume that $g^k_{l+1}(P) \neq 0$. 
	 We shall make elementary matrix operations on $M_P$ to calculate the rank of $M_P$. 
	 
	 For $\lambda = l+2, \ldots, c$, let 
	 \begin{multline*}
	 Z_{\lambda}(P):= \{ Q \in (f_1 = \cdots = f_l=0) \cap \Pi^0  : \\
	  g^k_{l+1}(P) g^i_{\lambda}(Q) - g^k_{\lambda}(P) g^i_{l+1}(Q) =0 \ \ (i= k+1, \ldots, n) \}. 
	 \end{multline*}
	Note that the first row $M_P^1$ and the $(\lambda-l)$-th row $M_P^{\lambda-l}$ of $M_P$ are linearly dependent 
	 if and only if $P \in Z_{\lambda}(P)$. 
	 By condition (Q2) for $J$ with $|J| =2$, there are at least $k-l$ nonzero elements of 
	 $G_{\lambda}(P):= \{ g^k_{l+1}(P) g^i_{\lambda} - g^k_{\lambda}(P) g^i_{l+1} : i= k+1, \ldots,n \}$ and 
	 they define $k-l$ free linear systems on $\Pi^0$.   
	 	Hence we obtain $Z_{\lambda}(P) = \emptyset$ and the two rows $M_P^1$ and $M_P^{\lambda -l}$ are linearly independent. 
		Thus, by elementary operations on $M_P$, we obtain a matrix of the following form; 
		\[
		\begin{pmatrix}
		g_{l+1}^k & \cdots & \cdots & g_{l+1}^n \\ 
		0 & h^{k+1}_{l+2} & \cdots & h^n_{l+2} \\ 
		\vdots & \vdots & & \vdots \\
		0 & h^{k+1}_{c} & \cdots & h^n_{c} 
		\end{pmatrix}(P).  
		\]
		By column exchange operations, we may assume that $h^{k+1}_{l+2}(P) \neq 0$ and repeat the process 
		 to 
		 \[
		 M'_P:= \begin{pmatrix}
		
		 h^{k+1}_{l+2} & \cdots & h^n_{l+2} \\ 
		  \vdots & & \vdots \\
		 h^{k+1}_{c} & \cdots & h^n_{c} 
		\end{pmatrix}(P).  
		 \]

		 Let $G'_{\lambda}(P):= \{ h^{k+1}_{l+2}(P) h^i_{\lambda} - h^{k+1}_{\lambda}(P) h^i_{l+2} : i= k+2, \ldots,n \}$. By condition (Q2) for $J$ with $|J|=3$, 
		  there are at least $k-l$ nonzero elements of $G'_{\lambda}(P)$ and 
		 they define free linear systems on $\Pi^0$. 
		By this, we again see that the first row and another row of $M'_P$ are linearly independent.  
		 
		 After repeating these elementary operations, we obtain a matrix of the form 
		 \[
		 \begin{pmatrix}
		 \alpha_{l+1} & &   & \\
		  & \ddots & \text{\huge *} & \text{\huge *} \\
		  \text{\huge 0}& & \alpha_{c} &
		 \end{pmatrix}
		 \]
		 for some $\alpha_{l+1}, \ldots, \alpha_c \in \mathbb{C}\setminus \{0 \}$ 
	and see that the rank of $M_P$ is $c-l$. Thus $C_X^*$ is nonsingular at $P \in \Pi^0$. 
	\vspace{5mm}
	
	Suppose that conditions (Q1) and (Q2) do not hold for some $I$. 
	We shall show that $X$ is not quasi-smooth. 
	We may again assume that $I = \{0, \ldots, k-1 \}$ and 
	$\Pi = (x_k = \cdots = x_n =0)$. 
	Moreover, since (Q1) and (Q2) do not hold, we may assume that, for some $l < \rho_I$, we have  
	$\Pi \not\subset (f_\lambda=0) $ for $\lambda = 1, \ldots, l$ and 
	$\Pi \subset (f_\lambda =0)$ for $\lambda = l+1, \ldots, c$. 
	Then the singular locus of $C^*_X$ on $\Pi^0$ can be described as 
	\[
	Z:= \{P \in (f_1 = \cdots = f_l=0) \cap \Pi^0 : \rk M_P < c-l  \}, 
	\]
	where $M_P$ is the matrix defined in (\ref{eq:M_P}). 
	By the hypothesis, we may also assume that there exists $J \subset \{l+1, \ldots, c \}$ 
	such that there are at most $k-l +|J| -2$ nonzero elements among 
	$\{g_{\lambda}^i : \lambda \in J, \ i=k, \ldots, n \}$. 
	This implies that there are at most $k-l+|J| -2$ nonzero columns of the matrix 
	$M_P^J:= \left( g_{\lambda}^i(P) \right)_{\substack{k \le i \le n \\ \lambda \in J}}$. 
	We can choose $J$ so that the number $|J|$ is minimal among such subsets of $\{l+1, \ldots, c \}$.  
	Then, by elementary operations as in the first part of the proof, 
	we can transfer $M_P^J$ to the form 
	\[
	\begin{pmatrix}
		h_{l+1}^k & \cdots  & \cdots & h_{l+1}^{k+|J| } & \cdots & h_{l+1}^n \\ 
		 & \ddots &  & \vdots & & \vdots \\ 	
		\text{\huge 0} & & h^{k+|J|-1}_{l+ |J|} & h^{k+|J| }_{l+ |J|} & \cdots & h^n_{l+|J|} 
		\end{pmatrix}(P).  
	\]
	Note that, on the bottom row, we have at most $k-l-1$ nonzero entries. 
	Hence we obtain 
	\begin{multline*}
	\dim (f_1 = \cdots = f_l =0) \cap (h^{k+|J|-1}_{l+|J|} = \cdots = h_{l+|J|}^{n} =0) \cap \Pi^0 \\
	\ge k-l-(k-l-1) = 1. 
	\end{multline*}
	Since the rank of $M_P^J$ is not maximal on the subset $(h^{k+|J|-1}_{l+|J|} = \cdots = h_{l+|J|}^{n} =0)$, we see that $C_X^*$ is singular along the above positive dimensional subset. 
	Hence $X$ is not quasi-smooth in this case. 
	This concludes the proof of Proposition \ref{prop:qsmcrit}.

	
	\end{proof}

In the following example, we use Proposition \ref{prop:qsmcrit} to check quasi-smoothness of a 
given WCI. 

\begin{example}
Let $X_{8,8,8} \subset \mathbb{P}(2^{(4)}, 3^{(5)}, 5^{(3)})$ be a general WCI of codimension $3$. 
We can check the quasi-smoothness of $X_{8,8,8}$ by Proposition \ref{prop:qsmcrit} as follows. 
Consider $I = \{4,5,6,7,8 \}$, that is, $a_4= \cdots=a_8 =3$. 
Then (Q1) does not hold for this $I$ and we have $k=5, l=0$ in (Q2). 
We can choose $\{e_{\mu, j} : j=1,2,3, \mu=1, \ldots, 5 \} \subset \{0,1,2,3, 9,10,11 \}$ 
so that (Q2) is satisfied for this $I$. 
We can similarly check that (Q2) holds for $I= \{9,10,11 \}$. 
For other $I$, we have (Q1), thus we see the quasi-smoothness of $X_{8,8,8}$. 

On the other hand, we see that $X'_{8,8,8} \subset \mathbb{P}(2^{(3)}, 3^{(4)}, 5^{(3)})$ is not quasi-smooth. 
Indeed, for $I= \{7,8,9 \}$, that is, $a_7=a_8=a_9=5$, neither (Q1) nor (Q2) hold. 
\end{example}


The following proposition treats the special situation where 
some weight of $\mathbb{P}$ divides none of the degrees of a WCI. 

\begin{proposition}\label{prop:fundanonempty}
	Let $X= X_{d_1, \ldots, d_c} \subset \mathbb{P}(a_0, \ldots, a_n)$ be a (well-formed)
	quasi-smooth WCI which is not a linear cone. 
	Assume that there exists $i_0$ such that $a_{i_0}$ does not divide $d_j$ for all $j$. 
	Let $H= \mathcal{O}_X(h)$ be the fundamental divisor on $X$, that is, an ample Cartier divisor on $X$ which generates $\Pic X$. Then 
	\begin{enumerate}
	\item[(i)] $X$ has a quotient singularity of type $1/a_{i_0}(c_1, \ldots, c_{n-c})$ for some 
	$c_1, \ldots, c_{n-c} \in \mathbb{Z}_{\ge 0}$ such that $\gcd(a_{i_0}, c_1, \ldots, c_{n-c}) =1$; 
	\item[(ii)] $a_{i_0} \mid h$. 
	As a consequence, we have $|H| \neq \emptyset$. 
	\end{enumerate}
	\end{proposition}
 
 \begin{proof}
	Let $f_1, \ldots, f_c \in \mathbb{C}[x_0, \ldots, x_n]$ be the defining equations of $X$ such that 
	$\deg f_j = d_j$ for $1 \le j \le c$ and $X = (f_1 = \cdots = f_c =0) \subset \mathbb{P}(a_0, \ldots, a_n)$, 
	where $\deg x_i = a_i$ for $0 \le i \le n$. 
 	By applying Proposition \ref{prop:qsmcrit} to $I = \{i_0 \}$, we see that there exist distinct integers $e_1, \ldots, e_c \in \{0, \ldots, \hat{i_0}, \ldots, n \}$ and positive integers
	$k_1 \ldots, k_c$ such that 
	$d_j = k_j a_{i_0} + a_{e_j}$ for $1 \le j \le c$, i.e.\ we can write
 	 \[
 	 f_j = x_{i_0}^{k_j} x_{e_j} + g_j 
 	 \]
 	 for $1 \le j \le c$, where $g_j$ is a weighted homogeneous polynomial of degree $d_j$. 
 	 
	By the inverse function theorem, we see that $X$ has a quotient singularity of type $1/a_{i_0}(a_0, \ldots, \hat{a_{i_0}}, \ldots,  \hat{a_{e_1}}, \ldots, \hat{a_{e_c}}, \ldots, a_n)$ at 
	  $P_{i_0}:= [0: \cdots:1: \cdots:0]$. 
	  We shall show that $g :=\gcd(a_0, \ldots, \hat{a_{e_1}}, \ldots, \hat{a_{e_c}}, \ldots, a_n) =1$. 	 Suppose that $g >1$. 
	 
	\begin{claim}\label{claim:c'}
	Up to a permutation on $\{1,\ldots,c\}$, we may choose $0 \le c' \le c$ with the following properties: 
	 \begin{enumerate}
	 \item[(*)] For $j= 1, \ldots, c'$, some monomial in $g_j$ does not contain 
	  any element of $\{x_{e_j}, \ldots, x_{e_c}  \}$. 
	  \item[(**)] For $j=c'+1, \ldots, c$, every monomial 
	  in $g_j$ contain some of $\{x_{e_{c'+1}}, \ldots, x_{e_c} \}$. 
	  \end{enumerate}
	 \end{claim}
	
	\begin{proof}[Proof of the claim]
	If (**) holds for all $j=1, \ldots, c$ and $\{x_{e_{1}}, \ldots, x_{e_c} \}$, then we put $c':=0$. 
	Otherwise there is some $j$ such that $1 \le j \le c$ and (*) holds for  $\{x_{e_{1}}, \ldots, x_{e_c} \}$. 
We then exchange $(f_1, e_1)$ and $(f_j, e_j)$ and  repeat the same process starting from $j=2$ till we obtain the claim, that is, check whether (**) holds for new $\{f_2, \ldots, f_c \}$ and $\{e_2, \ldots, e_c \}$ and so on. 
	\end{proof}

	Hence, for $1 \le j \le c'$, there exists a monomial in $g_j$ of the form $h_j = \prod_{i \neq e_j, \ldots, e_c} x_i^{b_i}$. 
	  Then we have $a_{e_j} \equiv \sum_{i \neq e_j, \ldots, e_c} b_i a_i \mod a_{i_0}$. 
	Thus we can check one by one that 
	\begin{equation}\label{eq:dividerel}
		g \mid a_{e_j} \text{ for }  1 \le j \le c'. 
	\end{equation}
	
	  
	  
	 
	  Now let $\Pi := (x_{e_{c'+1}} = \cdots = x_{e_c} =0) \subset \mathbb{P}$. We have $\Pi \subset \Sing \mathbb{P}$, in particular $\Pi \ne \P$.
	   
	  We also have $f_j |_{\Pi} \equiv 0$ for $c'+1 \le j \le c$ by the property (**) of $c'$. 
	  Thus we obtain 
	  \[
	  \dim \Pi \cap X \ge \dim \Pi - c' = \dim \mathbb{P} - c. 
	  \]
	  This contradicts the fact that $X \not\subset \Pi$ since $X$ is not a linear cone.  
	  Hence we obtain $g=1$, concluding the proof of Proposition \ref{prop:fundanonempty}.  
	   

	 \end{proof}

 The following proposition is useful for calculating the fundamental divisor of a WCI and is the motivation of the definition of \emph{$h$-regular pair} (see Definition \ref{def:h-regular}). 
 
 \begin{proposition}\label{prop:i1i2}
 	Let $X= X_{d_1, \ldots, d_c} \subset \mathbb{P}(a_0, \ldots, a_n)$ be a
 	quasi-smooth well-formed WCI which is not a linear cone. 
	Let $H= \mathcal{O}_X(h)$ be the fundamental divisor of $X$. 
	Assume that there exists $I= \{i_1, \ldots,  i_k \}$ such that $a_I:= \gcd(a_{i_1}, \ldots, a_{i_k}) >1$. 
	
	Then one of the following holds: 
	
	\begin{enumerate}
		\item[(i)] There exist distinct integers $p_1, \ldots, p_k$ such that $a_I \mid d_{p_1}, \ldots, d_{p_k}$. 
		\item[(ii)] $a_I \mid h$. 
		\end{enumerate} 
	 
	  \end{proposition}
 
 \begin{proof}
	 We apply Proposition \ref{prop:qsmcrit} to $I= \{i_1, \ldots, i_k \}$. 
	 Let \[
	 P_I := (x_0 = \cdots \hat{x}_{i_1} = \cdots = \hat{x}_{i_k} = \cdots = x_n =0) \subset \mathbb{P}
	 \]
	 be the $(k-1)$-dimensional stratum corresponding to $I$ and $P_I^0 \subset P_I$ be the open toric  stratum. 
	 
	 Suppose that condition (Q1) in Proposition \ref{prop:qsmcrit} holds, that is, there exist distinct integers 
	 $p_1, \ldots, p_k$ and non-negative integers $k_{j,i}$ for $j=1,\ldots, k$ and $ i \in I$ such that 
	$d_{p_j} = \sum_{i \in I} k_{j,i} a_i$. Then we have (i) in this case. 
	 
	 Suppose that (Q2) holds. Then there exist a permutation $p_1, \ldots, p_c$ of $\{1, \ldots, c \}$, an integer $l < \rho:= \min \{c, k \}$, non-negative integers $k_{j,i}$ for $j=1,\ldots, c$ and $i \in I$, and
	 distinct integers $e_{l+1}, \ldots, e_c$, which satisfy the following;  
	 \begin{itemize}
		 \item for $j=1, \ldots, l$, we have $\sum_{i \in I} k_{j,i} a_i = d_{p_j}$, 
		 \item for $j=l+1, \ldots, c$, we have $a_{e_j} + \sum_{i \in I} k_{j,i} a_i = d_{p_j}$. 
		 \end{itemize}

	We may assume that $(f_{p_j}=0) \cap P_I^{0} \neq \emptyset$ since $X$ is irreducible and 
	the linear system $|\mathcal{O}_{\mathbb{P}}(d_{p_j})|$ 
	does not have a fixed component. 
	Hence, on $p \in X \cap P_I^0$, the variety $X$ is analytic locally isomorphic to a quotient singularity of type 
	\[
	\frac{1}{a_I}(a_0, \ldots, \hat{a}_{i_1}, \ldots, \hat{a}_{i_k}, \ldots, 
	\hat{a}_{e_{l+1}}, \ldots, \hat{a}_{e_{c}}, \ldots, a_n) \times \mathbb{C}^{k-l}.
	\] 
	Now the proof is reduced to the following claim. 
	\begin{claim} 
	We have 
	$g:=  \gcd(a_0, \ldots,  
	 \hat{a}_{e_{l+1}}, \ldots, \hat{a}_{e_{c}}, \ldots, a_n) =1$.
	\end{claim}
	
	\begin{proof}[Proof of Claim]
	 Suppose that $g >1$.  We shall have a similar contradiction 
	as in the proof of Proposition \ref{prop:fundanonempty}. 
	As in Claim \ref{claim:c'}, up to a permutation of $\{1, \ldots, c \}$,
	 we may choose $c'$ with $l+1 \le c' \le c$ with the following properties: 
	
	\begin{enumerate}  
	 \item[(*)] For $j= l+1, \ldots, c'$, some monomial in $g_j$ does not contain 
	  any element of $\{x_{e_j}, \ldots, x_{e_c}  \}$. 
	  \item[(**)] For $j=c'+1, \ldots, c$, every monomial 
	  in $g_j$ contain some of $\{x_{e_{c'+1}}, \ldots, x_{e_c} \}$. 
	\end{enumerate} 
	Let $\Pi := (x_{e_{c'+1}} = \cdots = x_{e_{c}} =0) \subset \mathbb{P}$. 
	Then, as in Proposition \ref{prop:i1i2}, we see that 
	$f_j |_{\Pi} \equiv 0$ for $j= c'+1, \ldots, c$ and 
	$\Pi \subset \Sing \mathbb{P}(a_0, \ldots, a_n)$ since 
	$g \mid a_{i}$ for $i \notin \{e_{c'+1}, \ldots, e_{c} \}$.  
	 Thus we have $\dim \Pi \cap X \ge \dim \mathbb{P} -c$ as before and 
	 it contradicts that $X \not\subset \Pi$ since $X$ is not a linear cone. 
	 Thus we obtain the claim. 
	\end{proof}
	
	 The sheaf $\mathcal{O}_X(1)$ induces a generator of the class group of a quotient singularity of the above type. 
	Since the class group is a cyclic group of order $a_I$ as in (\ref{eq:classgroup}), we see that $a_I \mid h$. 
	Thus we finished the proof of Proposition \ref{prop:i1i2}.   
	 \end{proof}

	The following corollary restricts Proposition \ref{prop:i1i2} to the smooth case. 
	
\begin{corollary}(\cite[Lemma 2.15]{PS16})\label{cor:conditionsmooth}
 	Let $X= X_{d_1, \ldots, d_c} \subset \mathbb{P}(a_0, \ldots, a_n)$ be a
 	smooth WCI. 
	Assume that there exists $I= \{i_1, \ldots,  i_k \}$ such that $a_I:= \gcd(a_{i_1}, \ldots, a_{i_k}) >1$. 
	
	Then there exist distinct integers $p_1, \ldots, p_k$ such that $a_I \mid d_{p_1}, \ldots, d_{p_k}$.
	\end{corollary}
	
	\begin{proof}
		Since $X$ is smooth, the fundamental divisor of $X$ is $\mathcal{O}_X(1)$, that is $h =1$ 
		in the notation of Proposition \ref{prop:i1i2}. 
		Thus the statement follows from Proposition \ref{prop:i1i2}. 
		\end{proof}

\section{Regular pairs and Frobenius coin problem}


The following definition is motivated by Proposition \ref{prop:i1i2} and Corollary \ref{cor:conditionsmooth}.

\begin{definition}\label{def:h-regular}
Let $c \in \N$ and $n \in \Z_{\ge-1}$ be integers and  $(d; a)$ be a pair, where $d=(d_1, \ldots,d_c) \in \N_+^c$ and $a=(a_0, \ldots, a_n) \in \N_+^{n+1}$. Let $\bar{c}^+:= \{1, \ldots, c \}$ and $\bar{n}:= \{0, \ldots, n \}$. 

We say that $(d;a)$ is \emph{h-regular} for a positive integer $h$ if, for any subset $I= \{i_1, \ldots,  i_k \} \subset \bar{n}$ such that $a_I:=\gcd(a_{i_1}, \ldots, a_{i_k})>1$, one of the following holds:

 	\begin{enumerate}
		\item[(i)] There exist distinct integers $p_1, \ldots, p_k \in \bar{c}^+$ such that $a_I \mid d_{p_1}, \ldots, d_{p_k}$. 
		\item[(ii)] $a_I \mid h$. 
		\end{enumerate} 
	 
If a pair is $h$-regular for $h = 1$, we simply call it \emph{regular}.
\end{definition}

\begin{remark}
For technical reasons, in Definition \ref{def:h-regular} we admit the cases $c = 0$ or $n = -1$, i.e. pairs of the form $(d; \emptyset)$, $(\emptyset;a)$ and $(\emptyset,\emptyset)$.
\end{remark}

We need to fix some notation. If $(d;a)$ is a pair with $d=(d_1, \ldots,d_c) \in \N_+^c$ and $a=(a_0, \ldots, a_n) \in \N_+^{n+1}$, then we define
$$
\delta(d;a):=\sum_{j=1}^c d_j -\sum_{i=0}^n a_i.
$$
In the case where the pair $(d;a)$ comes from a well-formed quasi-smooth WCI $X=X_{d_1,\ldots,d_c} \subset \P(a_0,\ldots,a_n)$, we have $\omega_X \cong \m O_X(\delta(d;a))$.

Let $q$ be a prime number. Set $I_q:=\{i \in \bar{n} \ : \ q \mid a_i \}$ and $J_q:=\{ j \in \bar{c}^+  \ : \ q \mid d_j  \}$. We consider two new pairs obtained from $(d;a)$. 
The pair $(d^q;a^q)$ is given by 
$$
d^q :=((d_j /q)_{j \in J_q},(d_j)_{j \in \bar{c}^+ \setminus J_q}), \ a^q :=((a_i/q)_{i \in I_q}, (a_i)_{i \in \bar n \setminus I_q})
$$ 
in which we divided by $q$ all the divisible $d_j$ and $a_i$ and the pair $(d(q),a(q))$ is given by 
$$
d(q):=(d_j)_{j \in J_q}, \ a(q):=(a_i)_{i \in I_q} 
$$ 
in which only the divisible $d_j$ and $a_i$ appear. 
Note that
$$
\delta(d;a)= \delta(d^q ; a^q)+\frac{q-1}{q}\delta(d(q); a(q)).
$$

\begin{definition}\label{defn:cancellation}
For a pair $(d;a)$, we may uniquely choose  subsets 
$J_{(d;a)} = \{j_1, \ldots, j_l \} \subset \bar{c}^+$ and $I_{(d;a)} = \{i_1, \ldots, i_l \} \subset \bar{n}$ 
for some $l \in \mathbb{N}$ so that $d_{j_k} = a_{i_k}$ for $k=1, \ldots, l$ 
and $d_j \neq a_i$ for all $j \in \bar{c}^+ \setminus J_{(d;a)}$ and $i \in \bar{n} \setminus I_{(d;a)}$. 
We define a pair $(\tilde{d}; \tilde{a})$ by 
\begin{equation}\label{eq:tildedtildea}
(\tilde{d}; \tilde{a}) :=
 \left( (d_j)_{j \in \bar{c}^+ \setminus J_{(d;a)}}; (a_i)_{i \in \bar{n} \setminus I_{(d;a)}} \right), 
\end{equation}
that is, we cancel the doubles $(d_j, a_i)$ with $d_j = a_i$. 
\end{definition}

\begin{lemma}\label{lem:cancel}
The pair $(\tilde{d}; \tilde{a})$ is $h$-regular if $(d;a)$ is $h$-regular. 
\end{lemma}

\begin{proof}
Let $I:= \{i_1, \ldots, i_k \} \subset \bar{n} \setminus I_{(d;a)}$ 
be a subset with $a_I >1$. 
Since $(d;a)$ is $h$-regular, either (i) holds for  some $\{p_1, \ldots, p_k \} \subset \bar{c}^+$ 
or (ii) holds. 
In the latter case, there is nothing to check. 
Thus we consider the former case and need to find 
$p'_1, \ldots, p'_k \in \bar{c}^+ \setminus J_{(d;a)}$ such that $a_I \mid d_{p'_j}$ 
for $j=1, \ldots, k$. Let 
\[
J':= \{ j \in J_{(d;a)} : a_I \mid d_j \}, \ \ 
I':=  \{ i \in I_{(d;a)} : a_I \mid a_i \}. 
\]
Then we have $|I'| = |J'| =: l'$. 
Let $I'':= I \cup I'$. By $a_I = a_{I''}$, there exist distinct integers $p_1, \ldots, p_{k+l'} \in \bar{c}^+$  such that 
$a_I \mid d_{p_j}$ for $j=1, \ldots, k+l'$. 
Then the set $\{p_1, \ldots, p_{k+l'} \} \setminus J'$ contains $k$ elements 
$p'_1, \ldots, p'_k \in \bar{c}^+ \setminus J_{(d;a)}$ such that 
$a_I \mid d_{p'_j}$ for $j=1, \ldots, k$. 
Thus (i) holds for $(\tilde{d}; \tilde{a})$ and $I$. 
Hence we see that $(\tilde{d}; \tilde{a})$ is $h$-regular.    
\end{proof}

The following straightforward lemmas show how $h$-regular pairs are very suitable for inductive arguments.

\begin{lemma}\label{lemma:regular}
Let $(d;a)$ be an $h$-regular pair and $q$ be a prime not dividing $h$. Then the pairs $(d^q;a^q)$ and $(d(q);a(q))$ are $h$-regular. 
Hence $(d(q)/q;a(q)/q)$ is also $h$-regular.
\end{lemma}
\begin{proof}
We write the details for the pair $(d^q;a^q)$. The proof for $(d(q);a(q))$ is easier.

Let $I=\{i_1, \ldots,  i_k \} \subset \bar n$ such that $a_{I}:=\gcd(a_{i_1}, \ldots, a_{i_k})>1$. By the $h$-regularity of $(d;a)$,  we have either condition (i) of Definition \ref{def:h-regular}, i.e.\ there exist distinct integers $p_1, \ldots, p_k$ such that
$$ 
a_{I} \mid d_{p_1}, \ldots, d_{p_k}
$$ 
or (ii), i.e.\ $a_{I} \mid h$.

If  $q \mid a_{I}$, then we have $q \mid a_{i_\ell}$ for all $i_\ell \in I$ and $a_{I} \nmid h$.  Thus we have (i) and 
$$
\gcd \left(\frac{a_{i_1}}{q}, \ldots, \frac{a_{i_k}}{q}\right)= \frac{a_{I}}{q}  \mid \frac{d_{p_1}}{q}, \ldots, \frac{d_{p_k}}{q}
$$
as we wanted.

If $q \nmid a_{I}$, then 
$$
\gcd((a_i /q)_{i \in I'_q}, (a_i)_{i \in (I \setminus I'_q)})=a_{I}
$$
where $I'_q:= \{i \in I :  q \mid a_i  \}$.  If $a_{I} \mid h$, then there is nothing to prove, so we can assume that $a_{I} \nmid h$ and that (i) holds. 
Since $q \nmid a_{I}$, we get that $a_{I}$ divides $(d^q)_{p_j}$ for $j=1, \ldots ,k$. This concludes the proof.

\end{proof}

\begin{lemma}\label{lemma:h-regular}
Let $(d;a)$ be an $h$-regular pair and $q$ be a prime dividing $h$. Then 
$(d^q; a^q)$ is $h/q$-regular and $(d(q); a(q))$ is $h$-regular, hence  $(d(q)/q; a(q)/q)$ is $h/q$-regular.
\end{lemma}

\begin{proof}
We give the proof for the pair $(d^q; a^q)$. 
Consider a set $I= \{i_1, \ldots,  i_k \} \subset \bar n$ such that $a_I:=\gcd(a_{i_1}, \ldots, a_{i_k})>1$ and 
let $a^q_{I} := \gcd (a^q_i)_{i \in I}$ be the g.c.d.\ of the $a_i$'s in $(d^q;a^q)$. 

Assume first that $\gcd(a_I,q) = 1$, so that $a^q_I = a_I$. 
If $a_I \mid h$, we obtain that $a^q_I \mid h/q$ and we are done. 
If $a_I \nmid h$, then there exist distinct integers $p_1, \ldots, p_k$ such that $a_I \mid d_{p_1}, \ldots, d_{p_k}$. 
We have $a^q_I \nmid h/q$ and the $d^q_{p_j}$'s  work. 

If $a_I = qt$ for some positive integer $t$, then $a^q_I= t$. 
If $qt \mid h$,  we have $ t \mid h/q$. 
If $ qt \nmid h$,  then there exist distinct integers $p_1, \ldots, p_k$ such that $a_I \mid d_{p_1}, \ldots, d_{p_k}$. 
For the same integers,  we have $a^q_I \mid d_{p_1}/q, \ldots, d_{p_k}/q$ so the first condition of $h/q$-regularity is satisfied and we are done.


\end{proof}

\subsection{The Frobenius coin problem}\label{sec:Frobenius}

In this subsection we want to enlighten some interesting connections among Ambro--Kawamata's conjecture, regular pairs and the Frobenius coin problem.


\begin{question}[Frobenius coin problem]\label{question:Frobenius}
Given positive integers $a_0, \ldots, a_n$ such that $\gcd(a_0, \ldots, a_n)=1$, find the largest integer $G=G(a_0,\ldots,a_n)$ so that there do not exist nonnegative integers $x_0, \ldots, x_n$ satisfying
$$
G=a_0x_0 + \ldots + a_nx_n.
$$

Such $G$ is called the Frobenius number of $a_0, \ldots, a_n$.
\end{question}


 For $n=1$, it is classically known that
$$
G(a_0,a_1)=a_0a_1-a_0 -a_1.
$$
For $n \ge 2$, the problem is considerably harder: precise methods have been developed to compute $G(a_0,a_1,a_2)$ and some algorithms and (lower and upper) bounds are known for the general case (see for instance \cite{Johnson60} and \cite{BS62}). 


By Lemma \ref{lem:sectiondim}, Ambro--Kawamata's conjecture for smooth WCI would follow from the following purely arithmetic statement, which we believe to be of independent interest.

\begin{conjecture}\label{conj:regular}
Let $(d;a)=(d_1,\ldots,d_c;a_0,\ldots,a_n) \in \N^c \times \N^{n+1}$ be a regular pair such that $a_i \ne 1$ and $d_j \ne a_i$ for any $i,j$. Assume $c \le n$ and $\gcd(a_0,\ldots,a_n)=1$. Then 
$$
\delta(d;a) \ge G(a_0,\ldots,a_n).
$$
\end{conjecture}

One of the best known lower bounds for $G$ is given in \cite{Brauer}. Let $a_0,\ldots,a_n$ be positive co-prime integers, set $g_j:=\gcd(a_0,\ldots,a_j)$ for $j=0,\ldots,n$ and consider 
$$
Br(a_0,\ldots,a_n):= \sum_{j=1}^n a_{j} \frac{g_{j-1}}{g_j} - \sum_{i=0}^n a_i.
$$
Brauer proved that $Br(a_0,\ldots,a_n) \ge G(a_0,\ldots,a_n)$. Set $d_j:=a_{j} \frac{g_{j-1}}{g_j}$ for $j=1,\ldots,n$. Then it is easy to check that $(d;a):=(d_1,\ldots,d_n; a_0,\ldots,a_n)$ is actually a regular pair.

On the other hand, it is not difficult to see that, considering big prime numbers $p$ and $q$, the pair $(pq,6p,6q;2p,3p,2q,3q)$ is regular, $\delta(d;a) \ge G(a_0,\ldots,a_n)$, but $\delta(d;a) < Br(a_0,\ldots,a_n)$. 

This shows that regular pairs can give better bounds for the Frobenius number with respect to the known ones. 
For this reason, it seems to be a challenge and interesting problem to study Conjecture \ref{conj:regular}.

\begin{remark}\label{rmk:codimension2}
It is not difficult to check that Conjecture \ref{conj:regular} is true for $c=1,2$, which implies that the non-vanishing holds for a smooth WCI of codimension 1 or 2. For simplicity, we omit the detail in codimension $2$ case.  
\end{remark}

For $c=1$, a stronger and more general result is given in Lemma \ref{lemma:hypinequality}, which is the key step to prove Theorem \ref{thm:hypersurfaces}.

\section{Proof of Theorem \ref{thm:non-vanishing}}

Theorem \ref{thm:non-vanishing} is the combination of Corollary \ref{cor:smooth} and Corollary \ref{cor:singular} below.

\subsection{Smooth case}


%
%

The pair $(d;a)$ in the following lemma does not come from a non-empty WCI. 
Nevertheless this lemma is important in the proof of Proposition \ref{prop:regular}. 

\begin{lemma}\label{lemma:qdiv}
Let $(d; a) \in \N_+^c \times \N_+^{n+1}$ be a regular pair such that $a_i \neq d_j$ for any $i,j$. Let $q$ be a prime number such that $q \mid a_i$ 
and $q \mid d_j$ for any $i,j$. Then
$$
\delta(d;a) \ge cq.
$$
Moreover, if the equality holds, then $c=n+1$.
\end{lemma}

\begin{proof}
	 Note that $c \ge n+1$ which does not occur for a non-empty WCI. 

Assume first that  $q$ is the only prime dividing the $a_i$'s, that is for any $i=0,\ldots,n$, we have $a_i = q^{\alpha_i}$ for some $\alpha_i \ge 1$.  
We can assume that  $a_0 \ge a_1 \ge \ldots \ge a_n$. 
We can also order the $d_j$'s in such a way that $v_q(d_s) \ge v_q(d_t)$ for any $s \le t$, where $v_q(d_j)= \max \{e \in \N  :  q^e \mid d_j \}$. Then we have $a_{i} \mid d_{i+1}$ for any $i=0,\ldots,n$ and so 
$$
\sum_{j=1}^c d_j -\sum_{i=0}^n a_i = \sum_{k=1}^{c-n-1} d_{n+ 1 + k} + \sum_{i=0}^n (d_{i+1} - a_i) \ge cq
$$
and the equality is possible only if $c=n+1$, $d_j=2q$ and $a_i=q$ for any $i,j$.

\smallskip

Assume now that $q \neq 2$ and that $q$ and $2$ are the only primes dividing the $a_i$'s, that is for any $i=0,\ldots,n$ we have $a_i = 2^{\alpha_i} q^{\beta_i}$ for some $\alpha_i \ge 0$ and $\beta_i \ge 1$ such that $\alpha_i >0$ for at least one $i$.   We proceed by induction on $t= \max_{0 \le i \le n} \{\beta_i \}$, the greatest power of $q$ dividing at least one $a_i$. 

Suppose $t=1$. We can assume that  $v_2(a_i) \ge v_2(a_j)$ and $v_2(d_i) \ge v_2(d_j)$ for any $i \le j$.
Then again $a_i | d_{i+1}$ for any $i=0,\ldots, n$ and we conclude as before.

Suppose $t \ge 2$. Let $I_{q^t} := \{i \in \bar{n} : q^t \mid a_i \}$ and $J_{q^t} := \{j \in \bar{c}^+ :  q^t \mid d_j \}$. 
We consider the following pairs: $(d';a')$, where $d' = ((d_j /q)_{j \in J_{q^t}},(d_j)_{j \in \bar{c}^+ \setminus J_{q^t}})$ and 
$a' = ((a_i /q)_{i \in I_{q^t}},(a_i)_{i \in \bar{n} \setminus I_{q^t}}))$ and $(d'';a'')$, where $d'' = (d_j /q)_{j \in J_{q^t}}$ and $a'' = (a_i /q)_{i \in I_{q^t}}$. 
It is straightforward to check as in Lemma \ref{lemma:regular} that $(d';a')$ and $(d''; a'')$ are regular. 
Consider the regular pair $(\tilde{d}'; \tilde{a}')$ constructed in (\ref{eq:tildedtildea})
 which satisfies $\tilde{d}'_j \neq \tilde{a}'_{i}$ for any $i \in \bar{n} \setminus I_{(d';a')}, j \in \bar{c}^+ \setminus J_{(d';a')}$, where $I_{(d';a')} \subset \bar{n}$ and $J_{(d';a')} \subset \bar{c}^+$ 
 are the subsets defined in Definition \ref{defn:cancellation}.    
Let 
\[
m: = |\{j \in J_{q^t} : d_j/q = a_i \ (\exists i \in \bar{n} \setminus I_{q^t}) \}|, 
\]
\[ 
\bar{m} := |\{i \in I_{q^t} : d_j = a_i/q \ (\exists j \in \bar{c}^+ \setminus J_{q^t}) \}|. 
\]
Note that $|I_{(d';a')}|= |J_{(d';a')}| \le m+\bar{m}$. 
Let $k:=|J_{q^t}|$. 
By induction on $t$, 
 we may assume that we have $\delta(d';a')=\delta(\tilde{d}'; \tilde{a}') \ge (c-m -\bar m)q$ and $\delta(d'';a'') \ge kq$. 
 We have that $|I_{q^t}| \le |J_{q^t}|$ since $(d;a)$ is regular.  
  Since $m \leq k$ and $\bar m \leq |I_{q^t}| \leq k$, we obtain 
 \begin{multline}
\delta(d;a) = \delta(d';a') + (q-1)\delta(d'';a'') 
 \ge (c-2k)q + (q-1) kq \\
  =  cq -2kq + kq^2 - kq = cq+ kq(q-3) \ge cq, 
\end{multline}
because $q \ge 3$.  The equality is possible only if we have it for both $(\tilde{d}'; \tilde{a}')$ and $(d'';a'')$. This implies by induction on $t$ that $c=n+1$ in this case. 

\smallskip

We now pass to the general case. 
For any prime $p$, different from $q$ and $2$, let $e_p:= \max \{ e \in \mathbb{N} : p^{e} | a_i \ (\exists i) \}$. 
The proof is by induction on $D = \sum_p e_p$, where the index varies over all prime numbers different from $q$ and $2$. 
The case $D=0$ has already been treated in the first part of the proof. 
So assume $D \geq 1$ and that the inequality holds up to $D-1$. 
Consider $(d^p;a^p)$ and let 
\[
m_p:= |\{j \in \bar{c}^+ : d_j/p = a_i \ (\exists i \in \bar{n} \setminus I_p) \}|,  
\]  
\[
\bar{m}_p :=  |\{i \in \bar{n} : d_j = a_i/p \ (\exists j \in \bar{c}^+ \setminus J_p) \}|. 
\] 
Let us again consider the pair $(\tilde{d}^p; \tilde{a}^p)$
as in Definition \ref{defn:cancellation} by removing subsets $J_{(d^p; a^p)} \subset \bar{c}^+$ and 
$I_{(d^p; a^p)} \subset \bar{n}$. 
Then this satisfies the hypothesis $(\tilde{d}^p)_j \neq (\tilde{a}^p)_i$ for any $i$ and $j$. 
We again have that $|J_{(d^p; a^p)}| \le m_p +\bar{m}_p$. 
By induction on $D$, we obtain $\delta(d^p;a^p) = \delta(\tilde{d}^p;\tilde{a}^p) \geq (c-m_p-\bar{m}_p)q$. 
Now consider the pair $(d(p)/p;a(p)/p)$. 
Again by induction on $D$, we obtain $\delta(d(p)/p;a(p)/p) \geq sq$, 
where $s:= |\{j \in \bar{c}^+   :  p |d_j \}|$. 
We see that $m_p \le s$ by the definition of $m_p$. 
Let $s':=  |\{i \in \bar{n}   : p | a_i \}|$. 
We see that $s' \le s$ by the regularity of $(d;a)$ and that $\bar{m}_p \le s'$ by the definition of 
$\bar{m}_p$. Thus we have $\bar{m}_p \le s$.  
By these inequalities and $p \ge 3$, we conclude that  
\begin{multline*}
\delta(d;a) =  \delta(d^p;a^p) + (p-1)\delta(d(p)/p;a(p)/p) \geq (c-m_p-\bar m_p)q + (p-1)sq  \\
= cq + psq - m_p q -\bar m_p q -sq \geq cq + psq-3sq \geq cq
\end{multline*}
as we wanted. Again, the equality is possible only if $c=n+1$.
\end{proof}

\vspace{3mm} 

By using Lemma \ref{lemma:qdiv}, we prove the following key proposition.

\begin{proposition}\label{prop:regular}
Let $(d; a) \in \N_+^c \times \N_+^{n+1}$ be a regular pair such that $a_i >1$ and $a_i \neq d_j$ for any $i,j$. Then the following holds. 
\begin{enumerate}
\item[(i)] We have 
\begin{equation}\label{eq:deltaineq}
\delta(d;a)\ge c.
\end{equation}
\item[(ii)] If $\gcd(a_0,\ldots,a_n)=1$, then the equality holds only if $(d;a)$ is of the form $(6^{(s)}, 1^{(c-s)}; 2^{(s)},3^{(s)})$ for some integer $s$. 
\end{enumerate} 
\end{proposition}

\begin{proof}

(i) The proof is by induction on $n$, the case $n =0$ being obvious. 
We can assume that no prime divides every $a_i$, otherwise we are in the case of Lemma \ref{lemma:qdiv}. 
In particular, we may assume that there is a prime $q \ne 2$ which divides some $a_i$. 
Let 
$$
m:=| \{ j \in \bar{c}^{+} : d_j/q = a_i \ (\exists i \in \bar{n} \setminus I_q) \}|, 
$$
$$
\bar m:=  |\{i \in \bar{n} : d_j = a_i/q \ (\exists j \in \bar{c}^+ \setminus J_q, d_j \neq 1)  \}|, 
$$ 
$$
 \ell := |\{i \in \bar{n}  :  a_i=q  \}|,  
 \ \ \ \ s:= |\{ j \in \bar{c}^+ :  q| d_j  \}| = |J_q|. 
$$ 
We note that $m \leq s$ by definition and $\ell + \bar m \leq s$ by the regularity. 

\vspace{2mm}

\noindent({\bf Case 1}) Suppose that $\ell + m + \bar m \ge 1$. Then the pair $(d^q;a^q)$ has some redundant $a_i$, in the sense that $a_i/q=1$, $d_j/q=a_i$ or $d_j=a_i/q$ for some $i,j$. 
That is, we consider a regular pair $(\tilde{d}^q, \tilde{a}^q)$ and, 
by removing all $\tilde{a}^q_i =1$, we obtain 
a new regular pair $(\hat{d}^q; \hat{a}^q) \in \mathbb{N}_+^{\hat{c}} \times \mathbb{N}_+^{\hat{n}+1}$ 
for some $\hat{c} \le c$ and $\hat{n} \le n$. 
Note that $\hat{n} < n$ by the hypothesis $\ell  +m + \bar{m} \ge 1$. 
Let $\ell_1:= |\{j \in \bar{c}^+ : d_j =1 \}|$ and $\ell':= \min \{\ell, \ell_1 \}$. 
Then we see that $\hat{c} \ge c-m-\bar{m} - \ell'$ by the construction of $(\tilde{d}^q; \tilde{a}^q)$. 
Since we have $|\{i \in \bar{n} \setminus I_{(d^q; a^q)} : \tilde{a}_i^q =1 \}| = \ell - \ell'$, 
we obtain, by induction on $n$, that
$$
\delta(d^q; a^q) = \delta(\hat{d}^q; \hat{a}^q) - (\ell - \ell') \ge \hat{c}-(\ell - \ell') \ge c-\ell-m -\bar m.
$$
By applying Lemma \ref{lemma:qdiv} to $(d(q),a(q))$, we obtain
$$
\delta(d(q); a(q)) \ge sq. 
$$
By these and $\ell + m + \bar{m} \le 2s$, we obtain 
\begin{align}
\delta(d;a) &= \delta(d^q; a^q) + \frac{q-1}{q} \delta (d(q); a(q)) \label{eq:deltaeqc} \\ 
&\geq c-\ell - m  -\bar m + (q-1)s \geq c + qs - 3s \geq c  \nonumber 
\end{align}
since $q \geq 3$.

\vspace{2mm}
\noindent({\bf Case 2}) Suppose now that $\ell + m + \bar m = 0$. 
Then the pair $(d^q; a^q)$ satisfies the assumptions of the proposition. 
We note that
\begin{align*}
\delta(d;a) = \delta(d^q;a^q)+ \frac{q-1}{q}\delta(d(q); a(q)) > \delta(d^q;a^q) 
\end{align*}
since we have $\delta(d(q); a(q))>0$ by Lemma \ref{lemma:qdiv}. 
So we can replace the pair $(d; a)$ with $(d^q; a^q)$ without changing the number $c$ of $j$'s and we can repeat the argument from the beginning of the proof (possibly changing the prime $q$) till either we end up in ({\bf Case 1}) or we reach the situation of Lemma \ref{lemma:qdiv}. 
In both cases, we are done and obtain (\ref{eq:deltaineq}).  

\bigskip 

\noindent(ii) We now study when the identity holds in the case $\gcd(a_0,\ldots,a_n)=1$. 

Note that the case $n=1$ is clear, being equivalent to ask $a_0a_1-a_0-a_1=1$.

Assume $n\ge 2$ and let  $q \ne 2$ be a prime number such that $q \mid a_i$ for some $i$. 
We shall follow the proof of the inequality. In particular, we look at
$$
\delta(d;a)=\delta(d^q;a^q) + \frac{q-1}{q}\delta(d(q); a(q)).
$$ 
With the same notation as above, we note that the equality can hold only if we are in ({\bf Case 1}) and,  by Lemma \ref{lemma:qdiv}, the number $|\{i \in \bar{n} : q \mid a_i \}|$ must be equal to $s = |J_q|$.  
Moreover we obtain $q=3$ by (\ref{eq:deltaeqc}). This implies that the only possible prime numbers that divide at least one $a_i$ are 2 and 3.

We must also have $m=s$ and $\ell + \bar m=s$.
By $m=s$, we see that any $d_j/3 \in \mathbb{N}$ must be equal to some $a_i$ which is not divisible by $3$. 
Hence we can write
$$
(d;a)=(3\cdot 2^{\beta_1},\ldots,3 \cdot 2^{\beta_s},2^{\beta_{s+1}}, \ldots, 2^{\beta_c};  2^{\beta_1},\ldots,2^{\beta_s},3 \cdot 2^{\alpha_{s+1}}, \ldots, 3 \cdot 2^{\alpha_{n+1}})
$$
for some non-negative integers $\alpha_i$ and $\beta_i$. Then
$$
\delta(d;a)= \delta(2^{\beta_{s+1}}, \ldots, 2^{\beta_c}; 2^{\alpha_{s+1}}, \ldots, 2^{\alpha_{n+1}}) + \frac{2}{3}\delta(3\cdot 2^{\beta_1}, \ldots, 3\cdot 2^{\beta_s}; 3\cdot 2^{\alpha_{s+1}}, \ldots, 3 \cdot 2^{\alpha_{n+1}}).
$$
Also note that $\ell+ \bar{m} =s$ implies that $n+1-s = |\{ i \in \bar{n} : 3 \mid a_i \}| = \ell + \bar{m} = s$, thus $n+1 = 2s$. 
 By the regularity of $(d(3);a(3))$ and the assumption $d_j \ne a_i$ for any $i,j$, to have the equality $\delta(d(3) ; a(3))=3s$ we need  $\beta_j=1$ for $j=1,\ldots, s$ and $\alpha_i=0$ for $i=s+1, \ldots, n+1$, which implies
$$
\delta(d;a)=\delta(2^{\beta_{s+1}}, \ldots, 2^{\beta_c};1, \ldots,1) + 2s= \sum_{i=1}^{c-s}2^{\beta_{s+i}} - (n+1-s)+2s,
$$
i.e.\ $c=\delta(d;a)=\sum_{i=1}^{c-s}2^{\beta_{s+i}} +s$.

Hence, we must have $\beta_j=0$ for $j=s+1,\ldots,c$, which finishes the proof. 
\end{proof}

As a corollary of Proposition \ref{prop:regular}, we obtain 
the non-emptyness of $|\mathcal{O}(1)|$ and the smoothness of its general member 
on a smooth Fano or Calabi-Yau WCI.

\begin{corollary}\label{cor:smooth}
	Let $X:= X_{d_1, \ldots, d_c} \subset \mathbb{P}(a_0, \ldots, a_n)$ be 
	a well-formed smooth Fano or Calabi-Yau WCI which is not a linear cone.  
	Let $c_1:= |\{i \in \bar{n} : a_i =1 \}|$. 
	Then the following holds. 
	\begin{enumerate} 
		\item[(i)] We have $c_1 \ge c$. Moreover the equality is possible only if $X$ is Calabi-Yau of type $X_{6,\ldots,6} \subset \P(1^{(c)}, 2^{(c)},3^{(c)})$. 
	\item[(ii)] The linear system $|\mathcal{O}_X(1)|$ is non-empty 
	and its general member $H$ is smooth. 
	\end{enumerate}
	\end{corollary}

\begin{proof}
	\noindent(i) We may assume that $a_0 \le \cdots \le a_n$. 
	Thus we have $a_0 = \cdots = a_{c_1-1} =1$.   
	Since $X$ is smooth, we see that $(d_1, \ldots, d_c; a_{c_1}, \ldots, a_n)$ is regular. 
	By this and Proposition \ref{prop:regular} (i), we obtain 
	\[
	\delta(d_1, \ldots, d_c; a_{c_1}, \ldots, a_n) \ge c.
	\] 
	By the assumptions, we have $ 0 \ge \delta (d;a) \ge c-c_1 $ and this implies the former statement. 

Let $(d;a)$ be a regular pair which satisfies $c_1 = c$. 
Let $(\hat{d}; \hat{a})$ be the regular pair obtained by removing all $a_i =1$. 
Then $(\hat{d}; \hat{a})$ satisfies the hypothesis of  Proposition \ref{prop:regular} (ii) 
since $(d;a)$ defines a smooth WCI. 
Hence, by Proposition \ref{prop:regular}(ii), 
we see that $(\hat{d}; \hat{a}) = (6^{(c)} ; 2^{(c)}, 3^{(c)})$ and 
$(d;a) = (6^{(c)} ; 2^{(c)}, 3^{(c)}, 1^{(c)})$. 
 
	\vspace{2mm}
	
	\noindent(ii) 
	By the latter part of Proposition \ref{prop:regular}, we can assume that $X$ is not of the form $X_{6,\ldots,6} \subset \P(2^{(c)},3^{(c)},1^{(c)})$, otherwise the conclusion is immediate. In particular, we may assume $c_1 \ge c+1$.
	
	By (i), we see that $|\mathcal{O}_X(1)| \neq \emptyset$. 
	Since $X$ is smooth and well-formed, 
	we have $\Sing \mathbb{P}(a_0, \ldots, a_n) \cap X = \emptyset$. 
	Thus $H \cap \Sing \mathbb{P}(a_0, \ldots, a_n) = \emptyset$. 
	Hence it is enough to check $H$ is quasi-smooth at $P:= \Pi(p)$, where 
	$p \in \Pi^{-1}((x_0 = \ldots = x_{c_1-1}=0) \cap X)$ and $\Pi \colon \mathbb{A}^{n+1} \setminus \{0 \} \rightarrow \mathbb{P}(a_0, \ldots, a_n)$ 
	is the quotient map. 
	
	Set $H_i := X \cap (x_i =0)$ for $i = 0, \ldots, c_1-1$. 
	We shall look at the Jacobi matrices of $X$ and $H_i \subset \mathbb{P}(a_0, \ldots,a_{i-1}, a_{i+1}, \ldots, a_n)$. 
	Let $f_1, \ldots, f_c$ be the defining equations of $X$ such that $\deg f_j = d_j$. 
	For $i=0, \ldots, n$, set 
	\[
	\bold{v}_i(p) := \begin{pmatrix}
	{\partial f_1}/{\partial x_i} \\ 
	\vdots \\ 
	 {\partial f_c}/{\partial x_i} 
	 \end{pmatrix}
	(p). 
	\]
	The Jacobi matrix $J_X(p)$ and $J_{H_i}(p)$ of $X$ and $H_i$ can be written as  
	\[
	J_X(p) = \left( \bold{v}_0(p), \ldots, \bold{v}_n(p) \right), 
	\]
	\[
	J_{H_i}(p) = \left( \bold{v}_0(p), \ldots,\bold{v}_{i-1}(p), \bold{v}_{i+1}(p), \ldots, \bold{v}_n(p) \right) . 
	\] 
	Since $X$ is quasi-smooth, there exist linearly independent vectors 
 \[
 \bold{v}_{i_1}(p), \ldots, \bold{v}_{i_c}(p).
 \] 
 Since $c_1 \ge c+1$, we can choose $i$ so that $i \notin \{i_1, \ldots, i_c \}$. 
 Then we see that $H_i$ is quasi-smooth at $P:= \Pi(p)$. 
 Thus a general member $H$ is also quasi-smooth at $P$. 	
 	\end{proof}

\begin{remark}
Let $X \subset \mathbb{P}(a_0, \ldots, a_n)$ be a smooth WCI as in Corollary \ref{cor:smooth}.  
For $I \subset \bar{n}$ such that $a_I =1$, it may a priori happen that (Q1) does not hold, but (Q2) holds. 
That is why we make an argument as in Corollary \ref{cor:smooth} (ii). 
\end{remark}

\begin{remark}\label{rem:baselocus} 
Let $X_{d_1, \ldots, d_c}$ be a smooth WCI as in Corollary \ref{cor:smooth}. 
Motivated by a question by Andreas H\"{o}ring, we consider the description of the base locus $\Bsl |\mathcal{O}_X(1)|$. 

Up to reordering $d_1,\ldots,d_c$, we can assume that there is an integer $c' \le c$ with the following properties: 
 for $1 \le j \le c'$, there are weighted homogeneous polynomials $f_j(x_{c_1}, \ldots, x_n)$ of 
 degree $d_j$ and, for $c'+1 \le j \le c$, all monomials of degree $d_j$ contain one of the variables $x_0, \ldots, x_{c_1-1}$ 
 of weights $1$.  
 Since the base locus $\Bsl |\mathcal{O}_X(1)|$ is $(x_0 = \cdots = x_{c_1-1}=0) \cap X_{d_1, \ldots, d_c}$, 
 it is isomorphic to a general WCI $Y_{d_1, \ldots, d_{c'}} \subset \mathbb{P}(a_{c_1}, \ldots, a_n)$. 
 
 Thus the base locus is again a WCI. However this is not necessarily (quasi-)smooth in general. 
 We shall see this in Example \ref{ex:baselocus}. 
\end{remark}

\begin{example}\label{ex:baselocus}
Let $X:= X_{231, 231, 26} \subset \mathbb{P}:= \mathbb{P}(3,3,7,7,11,11,1^{(447)})$ be a general WCI. 
We can check that this is a smooth Fano WCI as follows: for $I = \{0,1 \}, \{2,3 \}$ or $\{4,5\}$, (that is, two variables of weights $3$, $7$ or $11$), 
we have (Q1) for $d_1 =231, d_2= 231$. 
Also, for $I= \{0,1,2,3 \}$ or $\{0,1,4,5 \}$, we have (Q1) for $d_1=231, d_2=231, d_3 = 26$ since $26 = 7 \cdot 2 + 3 \cdot 4 = 11 + 3 \cdot 5$. 
For $I = \{2,3,4,5 \}$, we have (Q2) for $d_1 = 231, d_2 = 231, d_3 =26 = 7 \cdot 2 + 11 + 1$. 
By Proposition \ref{prop:qsmcrit}, we see that $X$ is quasi-smooth, and smooth 
since $X \cap \Sing \mathbb{P} = \emptyset$.  

The base locus $\Bsl |\mathcal{O}_X(1)|$ is a WCI $Y:= Y_{231, 231, 26} \subset \mathbb{P}':= \mathbb{P}(3,3,7,7,11,11)$. 
This is not quasi-smooth. 
Indeed, for $I= \{2,3,4,5 \}$, neither (Q1) nor (Q2) holds because of the lack of suitable degree $26$ polynomials. 
In fact, $Y$ is a non-normal surface singular along a curve $(x_0=x_1=f_1=f_2=0) \subset \mathbb{P}'$, 
where $f_1, f_2$ are part of defining polynomials of degrees $231$ and $x_0, x_1$ are the variables of weights $3$.


Hence we can not expect smoothness of the base locus of the fundamental linear system even if 
it contains a smooth member. 
\end{example}

\begin{remark}
Let $W=W_{d_1,\ldots,d_c} \subset \P(a_0,\ldots,a_n)$ be a smooth WCI which is not a linear cone, where $a_i >1$ for any $i=0,\ldots,n$. By Corollary \ref{cor:smooth} we know that $W$ is not Fano. Then we can consider a WCI
$$
X=X_{d_1,\ldots,d_c} \subset \P(a_0,\ldots,a_n,1^{(\ell)})
$$
where $\ell=\delta(W) +1$. 
In this way $X$ is a smooth Fano with $-K_X=\mathcal O_X(1)$ and $\Bsl |\mathcal{O}_X(1)|$ is exactly $W$.
\end{remark}

In Corollary \ref{cor:smooth} we showed that for a smooth Fano WCI, the general member of the fundamental divisor is quasi-smooth. This is not true in general for a quasi-smooth Fano WCI as the following example shows.

\begin{example}\label{ex:funddiv}
Let $X=X_{35} \subset \P(5,7,2^{(k)},3^{(k)})$ where $k \ge 5$. Then $X$ is a quasi-smooth Fano WCI with fundamental divisor $\m O_X(6)$, but $X_{35,6} \subset \P(5,7,2^{(k)},3^{(k)})$ is not quasi-smooth. However, we see that a general member of $|\mathcal{O}_X(6)|$ 
has only terminal singularities. Indeed it has an isolated singularity at $[* : * : 0: \cdots :0]$ 
which is locally isomorphic to 
$0 \in (x_1^3 + \cdots + x_k^3 + x_{k+1}^2 + \cdots + x_{2k}^2 =0) \subset \mathbb{C}^{2k}$.  
\end{example}

It is also natural to look at the general element of $|-K_X|$ in the case of a Fano variety $X$. For instance, Shokurov and Reid (\cite{MR534602}, \cite{reid1983projective}) proved that a Fano $3$-fold with only canonical Gorenstein singularities admits an anticanonical member with only Du Val singularities.
Here we give an example of a smooth Fano WCI whose anticanonical members are singular (and not quasi-smooth). 
See also \cite[2.12]{MR2918165} for an example of a Fano $4$-fold with singular fundamental divisor.

\begin{example}(cf. \cite[Example 2.9]{MR3264677})
For $m \in \mathbb{Z}_{>0}$, 
let $X=X_{(2m+1)(2m+2)} \subset \mathbb{P}(1^{\left( 1+2m(2m+1) \right)}, 2m+1, 2m+2)$ 
be a weighted hypersurface of degree $(2m+1)(2m+2)$.  
Then we see that $-K_X = \mathcal{O}_X(2)$ and the linear system  $|{-i}K_X|$ does not contain 
	a smooth member for $i=1, \ldots, m$. 
	Thus we can not expect a smooth element of the pluri-anticanonical system 
	on a Fano manifold.  
	However, in the above example, we can find a member with only terminal singularities. 
	Moreover, the base locus of $|H|$ consists of a point. 
\end{example}

\begin{remark} It is well known that following the arguments in \cite[Section 5]{Ambro99} or \cite[Section 5]{Kawamata00} and assuming Conjecture \ref{conj:Ambro-Kawamata}, it is possible to show that the general element of $|-mK_X|$ has always only klt singularities for $m>0$ such that $-mK_X$ is Cartier (we thank Chen Jiang to have pointed this fact out to us).
\end{remark}

Finally, we also get the following corollary, which generalizes \cite[Corollary 4.2]{Przyjalkowski:2017ab} to any codimension. 

\begin{corollary}\label{cor:I(X)}
Let $X:= X_{d_1, \ldots, d_c} \subset \mathbb{P}(a_0, \ldots, a_n)$ be 
	a well-formed smooth Fano or Calabi-Yau WCI which is not a linear cone.  
	Let $c_1:= |\{i \in \bar{n} : a_i =1 \}|$. 
	Then $c_1 > I(X):= -\delta(d;a) = \sum_{i=0}^n a_i - \sum_{j=1}^c d_j.$
\end{corollary}
\begin{proof}
Consider the regular pair $(d;a)$ associated with $X$. 
We may assume that $a_0 \le \cdots \le a_n$, so that $a_0 = \cdots = a_{c_1-1} =1$. 
Let $(d';a')$ be the pair $(d;a_{c_1}, \ldots,a_n)$, where we took away the 1's from $(d;a)$. 
This pair is regular with no $a_i=1$ and so by Proposition \ref{prop:regular} we get
$$
\delta(d';a')\ge c >0,
$$
which implies
$$
\delta(d;a)=\delta(d';a')-c_1>-c_1,
$$
 i.e.\ $I(X) < c_1$, as we wanted.
\end{proof}

\subsection{General case}

The following is a key proposition to deduce the non-vanishing in the quasi-smooth Fano case. 

\begin{proposition}\label{prop:h-regular}
Let $h \in \mathbb{N}_+$ and $(d; a) \in \N_+^c \times \N_+^{n+1}$ be an $h$-regular pair with $c \ge 1$. If $a_i \nmid h$ for any $i=0,\ldots,n$ and $a_i \ne d_j$ for any $i,j$, then 
$$
\delta(d;a)> 0.
$$
\end{proposition}

\begin{proof}
Let us write $h=p_1^{\alpha_1}\cdots p_k^{\alpha_k}$, where the $p_i$'s are distinct prime numbers. The proof is by induction on $\alpha=\sum_{i=1}^k \alpha_i \ge 0$. 
If $\alpha=0$, then the pair $(d;a)$ is regular and the statement follows from Proposition \ref{prop:regular}, so we assume $\alpha \ge 1$.

Let $p$ be a prime number dividing $h$ and consider $(d^p; a^p)$. As usual,
$$
\delta(d;a)=\delta(d^p;a^p)+\frac{p-1}{p}\delta(d(p);a(p)).
$$

By Lemma \ref{lemma:h-regular}, $(d^p; a^p)$ and $(d(p)/p; a(p)/p)$ are $h/p$-regular. 
Note that there does not exist $i$ such that $a_i/p =1$ by the hypothesis $a_i \nmid h$. 
 Thus, after cancellation on $(d^p; a^p)$ (see Definition \ref{defn:cancellation}),  
we see that $(\tilde{d^p}; \tilde{a^p})$ and $(d(p)/p; a(p)/p)$ are $h/p$-regular and satisfy the hypothesis of the proposition.  

If $p\mid a_i$ or $p\mid d_j$ for some $i$ or $j$, then we obtain $\delta(d(p)/p; a(p)/p)>0$ 
by the induction hypothesis and conclude $\delta(d;a) >0$ by induction since we have either 
$\delta(d^p; a^p) = \delta(\tilde{d^p}; \tilde{a^p}) > 0$ 
or $(\tilde{d^p}; \tilde{a^p})$ is empty. 

If $p \nmid a_i$ and $p \nmid d_j$ for any $i,j$, then $\delta(d^p; a^p)>0$ by the induction hypothesis 
since $(d^p; a^p) = (d;a)$ is $h/p$-regular.  
Moreover $(d(p);a(p))$ is empty. Hence we can again conclude that $\delta(d; a)>0$.
\end{proof}

\begin{corollary}\label{cor:singular}
Let $X=X_{d_1,\ldots,d_c} \subset \P(a_0, \ldots, a_n)$ be a well-formed quasi-smooth WCI which is Fano or Calabi-Yau and which is not a linear cone. Then $|H| \neq \emptyset$ for any ample Cartier divisor $H$ on $X$.
\end{corollary}
\begin{proof}
Write $H=\m O_X(h)$. If there exists $i \in \bar{n}$ such that $a_i \mid h$, then we are done. 
Otherwise, we are in the situation of Proposition \ref{prop:h-regular}, and so the variety can not be Fano or Calabi-Yau 
since $(d; a)$ is $h$-regular by Proposition \ref{prop:i1i2}.
\end{proof}

\section{Weighted hypersurfaces}


The following lemma gives a proof of a generalized version of Conjecture \ref{conj:regular} in the case $c=1$.

\begin{lemma}\label{lemma:hypinequality}
Let  $a_0,...,a_n$ be positive integers, $n \ge 1$ and set 
 $$
 h :=\lcm_{i \ne j}(\gcd(a_i,a_j)).
 $$  
 
Assume that $a_i \nmid h$ for any $i$ and set $f :=\lcm(a_0,\ldots,a_n)$. Then 
$$
f-\sum_{i=0}^n a_i \ge \lcm(a_s,a_t) - a_s -a_t
$$
for any $s$ and $t$.
\end{lemma}

\begin{proof}
We first note that, for any proper subset $I$ of $\bar{n} := \{0, \ldots, n \} $, we have $f \ne \lcm_{i \in I}(a_i)$. In fact, suppose that the equality holds and let $k \in \bar{n} \setminus I$. For any prime power $p^e$ such that $e \ge 1$ and $p^e \mid a_k$, we have $p^e \mid f$. In particular, we have $p^e \mid a_\ell$ for some $\ell \in I$. This implies that $p^e \mid \gcd(a_k,a_\ell)$ and so $a_k \mid h$, which is a contradiction. In particular, $f \ge 2\lcm_{i \in I}(a_i)$.

The proof of the lemma is by induction on $n$ and the case $n=1$ is trivial, so assume $n\ge 2$. Let $s,t \in \bar{n}$ be such that $ s \ne t$. 
Then
$$
f - \sum_{i=0}^n a_i =\frac{f}{2} - a_s -a_t + \frac{f}{2} -\sum_{i \ne s,t} a_i \ge \lcm(a_s,a_t) - a_s -a_t + \lcm_{i \ne s,t} (a_i) -\sum_{i \ne s,t} a_i.
$$

If $n=2$, then we have 
$$
\lcm_{i \ne s,t} (a_i) -\sum_{i \ne s,t} a_i=0
$$
and we are done.

If $n \ge 3$, then we have  
$$
\lcm_{i \ne s,t} (a_i) -\sum_{i \ne s,t} a_i \ge \lcm(a_{s'},a_{t'}) -a_{s'} -a_{t'}
$$
for $s', t' \in \bar{n} \setminus \{s, t \}$ such that  $s' \ne t'$ by induction on $n$ and $\lcm(a_{s'},a_{t'}) -a_{s'} -a_{t'} \ge 0$ because we are assuming that $a_i \nmid h$ for any $i$.

\end{proof}

\begin{proposition}\label{prop:hypnonvanishing}
Let $X= X_d \subset \P=\P(a_0,\ldots,a_n)$ be a well-formed, quasi-smooth hypersurface of degree $d$ which is not a linear cone. Let $H$ be an ample Cartier divisor on $X$ such that $H - K_X$ is ample. 

Then $|H|$ is not empty.
\end{proposition}
\begin{proof}
Write $\m O_X(H)=\m O_X(h)$ for a positive integer $h$.



By Proposition \ref{prop:fundanonempty},  we can assume that $a_i \mid d$ for any $i$. 
Then $X$ is a Cartier divisor which intersects any stratum $P_{\{i, j \}}$ in some interior point. The condition of $H$ to be Cartier is then equivalent to 
$$
\lcm_{i \ne j} (\gcd(a_i,a_j)) \mid h. 
$$

If there exists $a_i$ such that $a_i \mid h$, then we are done. So assume  that $a_i \nmid h$ for any $i$ and let $f :=\lcm(a_0,\ldots, a_n)$. By Lemma \ref{lemma:hypinequality}, we get

$$
f-\sum_{i=0}^n a_i \ge \lcm(a_s,a_t) -a_s -a_t
$$
for any $s$ and $t$. Since $h > f-\sum_{i=0}^n a_i$ (because $H-K_X$ is ample and $f \mid d$) and $g:=\gcd(a_s,a_t) \mid h$ for any $s \ne t$, we can use the Frobenius number $G(a_s/g,a_t/g) = \frac{1}{g}(\lcm (a_s, a_t) - a_s - a_t)$ as in Section \ref{sec:Frobenius} to conclude that there are non-negative integers $\lambda_s, \lambda_t$ such that
$$
\lambda_s a_s + \lambda_t a_t=h,
$$
which implies that $|H|$ is not empty by Lemma \ref{lem:sectiondim}.
\end{proof}

\vspace{3mm}

In the following, we prove the base-point freeness on a Gorenstein weighted hypersurface. 

\begin{proposition}\label{prop:Fujita}
Let $X=X_d \subset \P=\P(a_0,\ldots,a_n)$ be a well-formed, quasi-smooth hypersurface of degree $d$ which is not a linear cone such that $K_X$ is Cartier. 
Let $H$ be the fundamental divisor of $X$ and $h$ be the positive integer such that 
$H =\m O_X(h)$.

Then $L=K_X+mH$ is globally generated for any $m \ge n$.
\end{proposition}

\begin{proof}

Suppose by contradiction that there is a point $p=[p_0: \cdots : p_n] \in \Bsl |L|$ and take $\ell$ such that $L=\m O_X(\ell)$.

Note that if $p_s \ne 0$ for some $s$, then $a_s \nmid h$, otherwise $x_s^e \in |L|$ for some positive integer $e$ and so $p \notin \Bsl |L|$. 
Also note that, for all $i \in \bar{n}$ such that $a_i \nmid h$, 
we have $a_i \mid d$ by Proposition \ref{prop:fundanonempty}. 

Assume first that there exists a unique $s \in \bar{n}$ such that $p_s \ne 0$. Since $p \in X$ and $a_s \nmid h$, we get that $a_s \mid d$. 
Let $f_d$ be the defining polynomial of $X_d$. If $f_d$ contains a monomial $x_s^{d/a_s}$, then we obtain $p \not\in X_d$ and this is a contradiction. 
If $f_d$ does not contain such a monomial, then it should contain a monomial of the form $x_s^{k} x_i$ for some $k>0$ and $i \neq s$ by the quasi-smoothness of $X_d$. 
Then we see that $a_s \mid a_i$ by $a_s \mid d$, and $X_d$ has a quotient singularity of index $a_s$. Thus we obtain $a_s \mid h$ and this is a contradiction. 

Hence we can assume that there exist $s$ and $t$  such that $s \ne t$, $p_{s} \ne 0$ and $p_{t} \ne 0$, thus $a_s, a_t \nmid h$. We have
$$
\ell= d -\sum_{i=0}^n a_i +mh = d -\sum_{a_i \nmid h} a_i - \sum_{a_i \mid h} a_i + mh.
$$

Assume that $- \sum_{a_i \mid h} a_i + mh \ge 1$. Since $a_i \mid d$ for all $i$ such that 
 $a_i \nmid h$,  we can apply Lemma \ref{lemma:hypinequality} to conclude that
$$
\ell > \lcm(a_s,a_t)-a_s-a_t,  
$$
which implies that $x_s^{e_s}x_t^{e_t} \in |L|$ for some non-negative integers $e_s$ and $e_t$. So we again have $p \notin \Bsl |L|$.

Assume now that $- \sum_{a_i \mid h}a_i + mh \le 0$. Then we can check that  $|\{i : a_i=h \}| \ge n-1$, because $m \ge n$.
Moreover, since $\P$ is well-formed, the greatest common factor of any $n$ weights is 1. By these, when $|\{i : a_i=h \}| =n$, we have $h=1$ and $\P=\P(a_0,1,\ldots,1)$ for some $a_0 >1$.  When $|\{i : a_i=h \}| = n-1$, we have 
$h=2$ and $\P=\P(1,1,2,\ldots,2)$. In both cases, we can check that $L$ is base point free, and we have derived a contradiction. 

\end{proof}

\section*{Acknowledgment} 

We thank D.\ Turchetti for many helpful conversations, C.\ Shramov for letting us know the recent preprint \cite{Przyjalkowski:2017ab}, 
and A.\ H\"{o}ring for asking us a question on the base locus of the fundamental linear system.
Part of this project has been realized while the second and third authors were visiting the Max Planck Institute for Mathematics and the Institute for Mathematical Science, NUS in Singapore. We are happy to thank them for their support during our stays. The second author was also partially supported by JST Tenure Track Program, JSPS KAKENHI Grant Numbers 
JP15J03158, JP16K17573, JP24224001.

%
%
%
%
%
%
%
%

\bibliographystyle{amsalpha}
\bibliography{Library}

\end{document}